\documentclass[11pt]{amsart}

\usepackage{amsmath,amssymb,amsthm,mathtools}
\usepackage{booktabs}
\usepackage{array}
\usepackage[margin=1.15in]{geometry}
\usepackage[colorlinks=true,linkcolor=blue!50!black,citecolor=blue!50!black,urlcolor=blue!50!black]{hyperref}
\usepackage{xcolor}
\usepackage[expansion=false]{microtype}
\usepackage{xurl}
\DeclareTextSymbolDefault{\textsection}{OMS}

\DeclareFontFamily{U}{wncy}{}
\DeclareFontShape{U}{wncy}{m}{n}{<->wncyr10}{}
\DeclareSymbolFont{mcy}{U}{wncy}{m}{n}
\DeclareMathSymbol{\Sha}{\mathord}{mcy}{"58}

\newcommand{\Q}{\mathbf{Q}}
\newcommand{\Z}{\mathbf{Z}}
\newcommand{\R}{\mathbf{R}}

\newcommand{\F}{\mathbf{F}}
\newcommand{\Qp}{\mathbf{Q}_p}
\newcommand{\Zp}{\mathbf{Z}_p}
\newcommand{\Qinf}{\Q_\infty}
\newcommand{\Sel}{\mathrm{Sel}}
\newcommand{\Gal}{\mathrm{Gal}}
\newcommand{\Hom}{\mathrm{Hom}}
\newcommand{\Reg}{\mathrm{Reg}}
\newcommand{\charid}{\mathrm{char}}
\newcommand{\ord}{\mathrm{ord}}
\newcommand{\cO}{\mathcal{O}}
\newcommand{\cE}{\mathcal{E}}
\newcommand{\cC}{\mathcal{C}}
\newcommand{\fP}{\mathfrak{P}}
\newcommand{\fq}{\mathfrak{q}}
\newcommand{\ff}{\mathfrak{f}}
\newcommand{\Dp}{D_p}
\newcommand{\msp}[1]{\left[#1\right]^{+}}
\newcommand{\Aplus}{A_{+}}

\theoremstyle{plain}
\newtheorem{theorem}{Theorem}[section]
\newtheorem{proposition}[theorem]{Proposition}
\newtheorem{lemma}[theorem]{Lemma}
\newtheorem{corollary}[theorem]{Corollary}
\newtheorem{conjecture}[theorem]{Conjecture}
\newtheorem*{theoremA}{Theorem A}
\newtheorem*{theoremB}{Theorem B}
\newtheorem*{observationC}{Observation C}
\theoremstyle{definition}
\newtheorem{definition}[theorem]{Definition}
\newtheorem{observation}[theorem]{Observation}
\theoremstyle{remark}
\newtheorem{remark}[theorem]{Remark}

\numberwithin{equation}{section}
\begin{document}

\title[$p$-primary Tate--Shafarevich groups of $y^2=x^3+39x$]{The $p$-primary Tate--Shafarevich groups of the rank-two CM elliptic curve $y^2=x^3+39x$ for primes $p<30000$}

\author{Xinrui Chang}
\author{Yuhang Xiao}

\date{November 27, 2026}


\begin{abstract}
Let $E$ be the elliptic curve $y^2=x^3+39x$, which has rank two and complex multiplication by $\Z[i]$. For each of the $1632$ primes $p\equiv 3\pmod 4$ with $7\le p<30000$, the good supersingular primes in this range, we prove with computer assistance that $\Sha(E/\Q)[p^\infty]=0$. The proof combines the signed main conjecture of Pollack and Rubin with an exact computation showing that the coefficient of $T^2$ in the plus $p$-adic $L$-function is a $p$-adic unit. With the published theorem of Coates, Liang and Sujatha at the ordinary primes $p\equiv1\pmod 4$ (supplemented at $p=17$ by a preprint of Banwait) and descent at $p=2,3,13$, this gives $\Sha(E/\Q)[p^\infty]=0$ for every prime $p<30000$. Only the supersingular part is new, and at the seven smallest supersingular primes the unit property can already be read off from Iwasawa invariants computed by Pollack and listed in the LMFDB. Numerically, this coefficient is congruent modulo $p$ to $-2d_p$ at all $1632$ primes, where $d_p$ comes from $p$-adic heights; this congruence is what the supersingular $p$-adic Birch--Swinnerton-Dyer conjecture predicts at the first $p$-adic digit. Pre-registered tests on two further curves, with predicted constants $8$ and $9/2$, are consistent with that prediction. We prove neither the finiteness of $\Sha(E/\Q)$ nor the rank-two Birch--Swinnerton-Dyer formula.
\end{abstract}

\maketitle

\tableofcontents

\section{Introduction}\label{sec:intro}

Let
\[
E\colon\ y^2=x^3+39x .
\]
This is a global minimal model; $E$ has conductor $N=48672=2^5\cdot3^2\cdot13^2$, discriminant $-2^6\cdot 3^3\cdot 13^3$, and complex multiplication by the maximal order $\Z[i]$ of $K=\Q(i)$. Its Mordell--Weil group is $E(\Q)\cong\Z^2\oplus\Z/2\Z$, generated by $P=(3,12)$, $Q=(27,144)$ and the $2$-torsion point $(0,0)$ (Section~\ref{sec:curve}). For curves of rank at least two neither the finiteness of the Tate--Shafarevich group nor the Birch and Swinnerton-Dyer (BSD) formula is known in a single case; what can be proved is the vanishing, or finiteness, of the $p$-primary part $\Sha(E/\Q)[p^\infty]$ for individual primes $p$. For five curves $y^2=x^3-Dx$ of rank two, among them $D=-39$, i.e.\ our curve $E$, Coates, Liang and Sujatha \cite{CLS1,CLS2} proved that $\Sha(E/\Q)[p^\infty]$ is finite at every good prime $p\equiv1\pmod4$ below $30000$, and zero except at four pairs $(D,p)$, one of which is $(-39,17)$ \cite[Thm.~1.3]{CLS2}; Banwait \cite[Prop.~7.4 and Cor.~7.5]{Banwait} has recently treated the four exceptional pairs in a preprint. For curves without complex multiplication, Stein and Wuthrich \cite[Thm.~1.1]{SteinWuthrich} proved $\Sha(E/\Q)[p]=0$ for $1{,}534{,}422$ pairs $(E,p)$ with $E$ of rank at least two and conductor at most $30000$ and $p$ a good ordinary prime, $5\le p<1000$, with surjective mod~$p$ representation. At supersingular primes, again for curves without complex multiplication, Perrin-Riou \cite{PR03} combined Kato's divisibility with numerical computation to determine the $p$-primary part of $\Sha$ for several curves, including curves of positive rank, and Stein and Wuthrich \cite[\S12.2]{SteinWuthrich} showed that $\Sha(E/\Q)[5]=0$ for the rank-two curve \texttt{1483a1}, at which $5$ is supersingular. For curves with complex multiplication, the theorem of Pollack and Rubin \cite{PollackRubin} provides an equality of characteristic ideals at supersingular primes in place of Kato's divisibility. The supersingular primes $p\equiv3\pmod4$ are not covered by the method of \cite{CLS1,CLS2}; Banwait \cite[\S14]{Banwait} lists this ``inert half'' as an open direction and points to the signed main conjecture of Pollack--Rubin \cite{PollackRubin}, the supersingular $p$-adic BSD conjecture, and work of Castella \cite{Castella25} as the relevant ingredients. In this paper we treat the supersingular primes below $30000$ for this one curve.

\subsection{Results}

\begin{theoremA}[computer-assisted]
Let $p\equiv 3\pmod 4$ be a prime with $7\le p<30000$; there are $1632$ such primes. Then
\[
\Sha(E/\Q)[p^\infty]=0\qquad\text{and}\qquad \Sha(E/\Q(i))[p^\infty]=0 .
\]
More precisely, for each such $p$ the Pontryagin dual $X^+$ of Kobayashi's plus Selmer group $\Sel_p^+(E/\Qinf)$ over the cyclotomic $\Zp$-extension is isomorphic to $\Zp^2$ with trivial action of $\Gal(\Qinf/\Q)$.
\end{theoremA}

The proof (Sections~\ref{sec:criterion} and~\ref{sec:computation}) consists of a criterion and a finite computation. The criterion (Proposition~\ref{prop:criterion}) says: if the coefficient $\Aplus$ of $T^2$ in Pollack's plus $p$-adic $L$-function $L_p^+(E,T)$ is a $p$-adic unit, then $\Sha(E/\Q)[p^\infty]=0$. It is deduced from published theorems of Pollack--Rubin \cite{PollackRubin} (the plus main conjecture, the absence of finite submodules, and an exact control theorem over $\Q$); our contribution is the short module-theoretic deduction, which we give in full. The computation determines $\Aplus\bmod p$ exactly, by a finite sum of modular symbols with denominators dividing $2$, for each of the $1632$ primes; it is nonzero in every case. The values are listed in the ancillary file (see ``Data and code'').

\begin{remark}[previously available data]\label{rem:lmfdb}
For the seven smallest of these primes, $p\in\{7,11,19,23,31,43,47\}$, the L-functions and Modular Forms Database (LMFDB) lists, on the page of this curve (label \texttt{48672.n2}, Cremona label \texttt{48672i1}), the analytic Iwasawa invariants of Pollack's plus and minus $p$-adic $L$-functions; according to the LMFDB these invariants were computed, rigorously, by R.~Pollack.\footnote{\url{https://www.lmfdb.org/EllipticCurve/Q/48672/n/2}; for the attribution and the statement on rigour see \url{https://www.lmfdb.org/EllipticCurve/Q/Source} and \url{https://www.lmfdb.org/EllipticCurve/Q/Reliability}; the definitions used are described at \url{https://www.lmfdb.org/knowledge/show/ec.iwasawa_invariants}. Pages consulted on 26 November 2026. The pages do not state explicitly which of the two listed values belongs to which sign; at the six primes other than $11$ the two values are equal, and at $p=11$ our computation ($[T^2]L_p^-\equiv0\pmod{11}$) identifies $\lambda^-=4$.} The plus invariants are $(\mu^+,\lambda^+)=(0,2)$ at all seven primes. Since $\ord_T L_p^+\ge2$ (Lemma~\ref{lem:order}), $\mu^+=0$ and $\lambda^+=2$ together mean that $\Aplus$ is a unit. Combined with Proposition~\ref{prop:criterion}, these publicly available data therefore already give Theorem~A for these seven primes, assuming that the normalisation of $L_p^+$ used there differs from ours by a $p$-adic unit. What is new here is the unit property, and hence the vanishing of $\Sha(E/\Q)[p^\infty]$, at the remaining $1625$ primes, and the explicit residues $\Aplus\bmod p$ in a fixed normalisation at all $1632$ primes; at the seven small primes our computation agrees with the LMFDB data. (At $p=11$ the LMFDB lists the two $\lambda$-values $2$ and $4$; the value $4$ for the minus sign is consistent with $[T^2]L_p^-\equiv0\pmod{11}$, which we also find, but the minus sign plays no role in this paper.)
\end{remark}

Combining Theorem~A with results from the literature and with descent arguments at the three bad primes gives:

\begin{theoremB}[combined; partly computer-assisted]
For every prime $p<30000$, $\Sha(E/\Q)[p^\infty]=0$. For the $1610$ good ordinary primes this is the theorem of Coates, Liang and Sujatha \cite[Thm.~1.3]{CLS2}, together with Banwait's preprint \cite[Prop.~7.4]{Banwait} at $p=17$; it is quoted, not reproved, here.
\end{theoremB}

Only one part of Theorem~B is new. The sources for the $\pi(30000)=3245=3+1632+1610$ primes are listed in Table~\ref{tab:sources}; the proof is in Section~\ref{sec:combined}. The theorem \cite[Thm.~1.3]{CLS2} used for the good ordinary primes names this curve; at $p=17$ it proves finiteness only, and the vanishing there rests on the preprint \cite{Banwait}. We have not repeated the computations behind these results.

\begin{table}[ht]
\centering
\small
\renewcommand{\arraystretch}{1.25}
\begin{tabular}{@{}>{\raggedright}p{2.6cm} r >{\raggedright}p{4.9cm} >{\raggedright\arraybackslash}p{3.6cm}@{}}
\toprule
primes & number & source & status \\
\midrule
$p=2$ & $1$ & complete $2$-isogeny descent (Prop.~\ref{prop:p2}) & standard descent; computer-assisted (mwrank) \\
$p=3$ & $1$ & $3$-descent over $\Q(i)$ via $K(E[3])$ and a certified class group (Prop.~\ref{prop:p3}) & descent; computer-assisted (PARI, unconditional) \\
$p=13$ & $1$ & Rubin's class number formula \cite[Thm.~3.3]{Rubin91}, elliptic units, descent (Prop.~\ref{prop:p13}) & published theorem plus descent; computer-assisted \\
$p\equiv1\ (4)$, $p\ne13,17$ & $1609$ & Coates--Liang--Sujatha \cite[Thm.~1.3]{CLS2} & published, computer-assisted; not repeated here \\
$p=17$ & $1$ & finiteness: \cite[Thm.~1.3]{CLS2}; vanishing: Banwait \cite[Prop.~7.4]{Banwait} & preprint \\
$p\equiv3\ (4)$, $7\le p<30000$ & $1632$ & Theorem~A (Sections~\ref{sec:criterion}--\ref{sec:computation}) & \textbf{this paper}; computer-assisted \\
\midrule
total & $3245$ & & \\
\bottomrule
\end{tabular}
\caption{Where the vanishing of $\Sha(E/\Q)[p^\infty]$ comes from, for $p<30000$.}\label{tab:sources}
\end{table}

\subsection{Method}
Let $p\equiv3\pmod4$, $p\ge7$. Then $p$ is inert in $K$, $E$ has good supersingular reduction at $p$ and $a_p=0$. Kobayashi \cite{Kobayashi} defined plus and minus Selmer groups over the cyclotomic $\Zp$-extension $\Qinf/\Q$ whose duals are $\Lambda$-torsion, and Pollack \cite{Pollack} defined plus and minus $p$-adic $L$-functions $L_p^\pm\in\Lambda=\Zp[[T]]$. For CM curves, Pollack and Rubin \cite[Thm.~7.3]{PollackRubin} proved the signed main conjecture $\charid_\Lambda(X^\pm)=(L_p^\pm)$, building on Rubin's main conjecture for imaginary quadratic fields \cite{Rubin91}. Using in addition that $X^+$ has no nonzero finite submodule and that $\Sel_p(E/\Q)\to\Sel_p^+(E/\Qinf)^\Gamma$ is an isomorphism \cite[Lemma~8.3]{PollackRubin}, we show that $T^2$ divides $L_p^+$ (because the rank is at least two) and that $X^+\cong\Zp^2$ as soon as $\Aplus=[T^2]L_p^+$ is a unit. The coefficient $\Aplus\bmod p$ is determined by the first layer of the cyclotomic tower: by the interpolation property of $L_p^+$ and Birch's formula it equals $-\sum_a \msp{a/p^2}\binom{k(a)}{2}$, a sum over $(\Z/p^2\Z)^\times$ of plus modular symbols weighted by a binomial coefficient of a discrete logarithm (Lemma~\ref{lem:finite-layer}). Using an exact integer table of the $104832$ Manin symbols of level $48672$ and Manin's continued-fraction algorithm, this costs $O(p^2\log p)$ integer operations per prime; the whole range took about $67$ minutes on one processor.

\subsection{A numerical observation}
The second part of the paper (Sections~\ref{sec:leading}--\ref{sec:control}) is numerical and is not used in the proof of Theorem~A or~B. Let $\Reg_p=R_1\,\omega+R_2\,F\omega$ be the $\Dp$-valued $p$-adic regulator of $P,Q$ (Perrin-Riou; we use the normalisation of Stein--Wuthrich \cite[Lemma~4.2]{SteinWuthrich}), written in the basis $(\omega,F\omega)$ of the Dieudonn\'e module, where $\omega=dx/2y$ and $F$ is crystalline Frobenius. In all computed cases $R_1\in p^2\Zp$, and we put $d_p=(R_1/p^2)\bmod p$.

\begin{observationC}[numerical evidence]
For all $1632$ primes $p\equiv3\pmod4$ with $7\le p<30000$,
\[
\Aplus\equiv-2\,d_p\pmod p .
\]
At $p=7$ and $p=11$ the corresponding congruence modulo $p^2$ also holds.
\end{observationC}

The constant $2$ is the value of $\prod_\ell c_\ell/\#E(\Q)_{\rm tors}^2=8/4$, and $\#\Sha(E/\Q)$ is expected to be $1$. We explain in Section~\ref{sec:leading} that the congruence $\Aplus\equiv-C\cdot\#\Sha\cdot d_p\pmod p$ is what the supersingular $p$-adic BSD conjecture of Bernardi and Perrin-Riou \cite{BPR93,PR03}, in the form stated by Stein and Wuthrich \cite[Conj.~5.2]{SteinWuthrich} (see also Sprung \cite{Sprung}), predicts at the first $p$-adic digit. The observation on $E$ was made in the data before this derivation was written down. To test it on curves where the predicted constant is not $2$, we fixed predictions in advance for the curves $y^2=x^3+25721x$ and $y^2=x^3+51137x$, whose complex BSD quotients $L''(E,1)/(2\Omega_E\Reg_\infty)$ are numerically $8$ and $9/2$, and then computed $\Aplus$ at the first seven, respectively six, good supersingular primes $p\ge7$. All observed values are consistent with the predicted constants $8$ and $9/2$ and not with the alternatives that were fixed at the same time (Section~\ref{sec:control}). On the ordinary side, Calonge Fabregat \cite{CalongeFabregat} has tested, for rank-two curves $y^2=x^3-Dx$ at split primes $p\le41$, the conjecture that an exact $p$-adic BSD quotient equals the complex quotient $L''(E,1)/(2\Omega_E\Reg_\infty)$: it holds modulo $p^3$ in all $5520$ pairs $(E,p)$ to which the method of \cite{CalongeFabregat} applies, coming from $1898$ curves with conductors up to $3.9\cdot10^{12}$ and including curves with analytic $\Sha$ of order $4$ and $9$ \cite[\S1.3 and Table~1]{CalongeFabregat}. Our test is the analogous one at inert primes, at much lower precision and on three curves.

\subsection{What this paper does not claim}\label{subsec:notclaim}
\begin{itemize}
\item We do not prove that $\Sha(E/\Q)$ is finite: Theorem~B concerns finitely many primes.
\item We do not prove the BSD formula for $E$, nor any statement about $L''(E,1)$. That $L''(E,1)/(2\Omega_E\Reg_\infty)$ equals $2$, so that the analytic order of $\Sha$ is $1$, is a numerical identification to $60$ digits, not a theorem.
\item The theorems of this paper are computer-assisted: they depend on the correctness of the programs and of the software libraries named in the text (eclib, PARI/GP, SageMath).
\item Part II (Observation~C and Section~\ref{sec:control}) is numerical evidence for a conjecture. It proves nothing, it agrees with the prediction only at the first $p$-adic digit (at the second digit in four cases), and the analytic orders $1$, $4$, $9$ used there are numerical identifications.
\item Only the supersingular primes in Theorem~B are new, and at the seven smallest of them the conclusion already follows from data in the LMFDB together with Proposition~\ref{prop:criterion} (Remark~\ref{rem:lmfdb}); the remaining primes are covered by published results, one preprint, and standard descent.
\end{itemize}

\subsection*{Organisation}
Section~\ref{sec:curve} fixes notation. Section~\ref{sec:criterion} proves the criterion, Section~\ref{sec:computation} describes the computation, and Sections~\ref{sec:bad} and~\ref{sec:split} treat the bad and the ordinary primes; Theorem~B is proved in Section~\ref{sec:combined}. Sections~\ref{sec:leading} and~\ref{sec:control} contain the numerical part, and Section~\ref{sec:limits} collects the limitations.

\section{The curve and notation}\label{sec:curve}

\subsection{Basic data}
The following facts were computed with PARI/GP \cite{PARI2} and SageMath \cite{sagemath}; they are standard and have been checked with both systems.

\begin{itemize}
\item The model $[a_1,a_2,a_3,a_4,a_6]=[0,0,0,39,0]$ is minimal, with discriminant $\Delta=-3796416=-2^6\cdot3^3\cdot13^3$ and conductor $N=48672=2^5\cdot3^2\cdot13^2$. Since $\Delta<0$, $E(\R)$ is connected.
\item The bad primes $2,3,13$ all have Kodaira type III and Tamagawa number $2$; the exponents of the conductor are $5,2,2$. Hence $\prod_\ell c_\ell=8$.
\item $E(\Q)_{\rm tors}=\{O,(0,0)\}\cong\Z/2\Z$; the global root number is $+1$.
\item $E(\Q)$ has rank $2$, with basis $P=(3,12)$, $Q=(27,144)$ of the free part (a complete $2$-descent with mwrank \cite{eclib} gives rank exactly $2$, and saturation shows that $P,Q$ generate). The N\'eron--Tate regulator is $\Reg_\infty=2.9178431649936463423\ldots$ (heights not divided by $2$).
\item The real period of the N\'eron differential $\omega=dx/2y$ is $\Omega_E=\int_{E(\R)}|\omega|=1.4838536397546222008\ldots$, the least positive real period.
\item Numerically, $L''(E,1)=17.318608802435875071\ldots$, and $L''(E,1)/(2\Omega_E\Reg_\infty)=2.000\ldots$ to $60$ digits.
\end{itemize}

$E$ has complex multiplication by $\Z[i]$, via $[i](x,y)=(-x,iy)$. The quadratic twist of $E$ by $K/\Q$ is $-y^2=x^3+39x$, which is isomorphic to $E$ over $\Q$ via $x\mapsto -x$. Hence $\operatorname{rank}E(K)=2\operatorname{rank}E(\Q)=4$. Over $K$, $E$ has reduction type $I_2^*$ (Tamagawa number $4$) at $1+i$ and type III (Tamagawa number $2$) at $3$ and at $3\pm2i$, and $E(K)_{\rm tors}=\Z/2\Z$.

For a prime $p\nmid 78$, $p$ is inert in $K$ if and only if $p\equiv3\pmod4$; then $E$ has good supersingular reduction at $p$ and $a_p=0$ (Deuring). We call these primes the \emph{inert} (or supersingular) primes; for $p\equiv1\pmod 4$, $p\ne13$, the reduction is good ordinary and $p$ splits in $K$. Below $30000$ there are $3245$ primes: the three bad primes $2,3,13$; $1632$ inert primes $7\le p\le29983$; and $1610$ split primes $5\le p\le29989$, $p\neq 13$.

\subsection{Iwasawa-theoretic notation}
Fix an inert prime $p\ge7$. Let $\Qinf/\Q$ be the cyclotomic $\Zp$\nobreakdash-extension, $\Gamma=\Gal(\Qinf/\Q)$, and let $\gamma\in\Gamma$ be the topological generator whose image under the cyclotomic character, projected to $1+p\Zp$, is $1+p$. Put $\Lambda=\Zp[[\Gamma]]\cong\Zp[[T]]$ with $T=\gamma-1$. For $a\in(\Z/p^2\Z)^\times$ write $\langle a\rangle=a/\omega_{\rm T}(a)\in1+p\Zp$, where $\omega_{\rm T}$ is the Teichm\"uller character, and define $j(a)\in\Z/p\Z$ by $\langle a\rangle\equiv(1+p)^{j(a)}\pmod{p^2}$.

$\Sel_p(E/F)$ denotes the $p^\infty$-Selmer group, and $\Sel_p^\pm(E/\Qinf)$ the plus and minus Selmer groups of Kobayashi \cite{Kobayashi}, in the normalisation of \cite[Def.~3.3]{PollackRubin}, and we write $X^\pm$ for the Pontryagin dual $\Hom(\Sel_p^\pm(E/\Qinf),\Qp/\Zp)$. Let $L_p^\pm=L_p^\pm(E,T)\in\Lambda$ be Pollack's signed $p$-adic $L$-functions \cite[Thm.~5.6]{Pollack}, normalised by the interpolation formulas \cite[(10)--(12), p.~460]{PollackRubin} with $\Omega_E$ as above. We write
\[
\Aplus=[T^2]L_p^+\in\Zp,\qquad\text{and also }\Aplus\bmod p .
\]
Replacing $\gamma$ by $\gamma^v$ ($v\in\Zp^\times$) multiplies $\Aplus$ by $v^{-2}$ once the constant and linear terms vanish (Lemma~\ref{lem:order}); so the residue depends on our choice of $\gamma$, but whether it is a unit does not.

\subsection{Modular symbols}
Let $f$ be the newform of level $48672$ attached to $E$. We use the normalisation
\[
2\pi i\int_{i\infty}^{r}f(z)\,dz=\msp{r}\,\Omega^+ +[r]^-\,\Omega^-,\qquad r\in\Q,
\]
with $\Omega^+=\Omega_E$ and $\Omega^-\in i\R$, so that $\msp{r}\in\Q$ and $\msp{0}=L(E,1)/\Omega_E$. This is the convention of SageMath's \texttt{modular\_symbol}; it refers to the period of $E$ itself, so the Manin constant does not enter.

\section{The criterion at supersingular primes}\label{sec:criterion}

Throughout this section $p\equiv3\pmod4$ and $p\ge7$.

\subsection{Published input}

\begin{theorem}[Pollack--Rubin]\label{thm:PR}
\begin{enumerate}
\item[(i)] {\rm\cite[Thm.~7.3, p.~461]{PollackRubin}} $X^\pm$ is a finitely generated torsion $\Lambda$-module and $\charid_\Lambda(X^\pm)=L_p^\pm\Lambda$.
\item[(ii)] {\rm\cite[Thms.~4.3 and~6.3, pp.~455--457; see the proof of Thm.~8.2, p.~463]{PollackRubin}} $X^\pm$ has no nonzero finite $\Lambda$-submodule.
\item[(iii)] {\rm\cite[Lemma~8.3, p.~462]{PollackRubin}} Restriction induces an isomorphism $\Sel_p(E/\Q)\xrightarrow{\ \sim\ }\Sel_p^\pm(E/\Qinf)^\Gamma$. Dually, $X^\pm/TX^\pm\cong\Hom(\Sel_p(E/\Q),\Qp/\Zp)$.
\end{enumerate}
\end{theorem}

The standing hypotheses of \cite{PollackRubin} (\S1, p.~450, repeated at the beginning of \S8, p.~461) are: $E/\Q$ has complex multiplication by the maximal order of an imaginary quadratic field $K$, and $p>2$ is a prime of good supersingular reduction. There are no conditions on the conductor, the Tamagawa numbers or the torsion. Both hypotheses hold for our $E$ and $p$.

\begin{remark}[on the use of Rubin's main conjecture]\label{rem:rubin}
The proof of \cite[Thm.~6.1, p.~457]{PollackRubin} invokes \cite[Thm.~4.1(ii)]{Rubin91}, twisted by $\rho^{-1}$, where $\rho$ is the character giving the action of $G_K$ on $E[p^\infty]$. For primes not split in $K$, \cite[Thm.~4.1(ii), p.~37]{Rubin91} gives in general only a divisibility of characteristic ideals; equality of the $\chi$-components is asserted for those irreducible $\Zp$-representations $\chi$ of $\Delta=\Gal(K_0/K)$ that are nontrivial on the decomposition group of $\mathfrak p$ in $\Delta$, where $\mathfrak p$ is the fixed prime of $K$ above $p$ (here $\mathfrak p=p\Z[i]$; see also \cite[p.~55 and Thm.~10.6, p.~57]{Rubin91}). In the situation of \cite{PollackRubin} this condition holds: there $K_0=K(E[p])$, $[K_0:K]=p^2-1$ is prime to $p$, and the Hilbert class field of $K$ is $K$ itself; $p$ is totally ramified in $K(E[p^\infty])/K$ \cite[\S1]{PollackRubin}, so the decomposition group of $\mathfrak p$ in $\Delta$ is all of $\Delta$; and the component used in \cite[\S6]{PollackRubin} is that of the irreducible $\Zp$-representation $\chi$ of $\Delta$ whose extension of scalars to $\cO_p=\Zp[i]$ contains $\rho|_\Delta\colon\Delta\cong(\Z[i]/p)^\times\to\Zp[i]^\times$. Since $\rho|_\Delta$ is faithful of order $p^2-1>1$, $\chi$ is nontrivial on $\Delta$, which is the decomposition group. Also, the standing hypothesis of \cite[\S1, p.~27]{Rubin91}, that $p$ does not divide the number of roots of unity in the Hilbert class field of $K$, holds for every odd $p$. This verification is not written out in \cite{PollackRubin}; we include it for completeness.
\end{remark}

\subsection{The order of vanishing}

\begin{lemma}\label{lem:order}
$\ord_T L_p^+\ge 2$. In particular $L_p^+(0)=0$ and the coefficient of $T$ in $L_p^+$ vanishes.
\end{lemma}

\begin{proof}
Since $\operatorname{rank}E(\Q)\ge2$, $\Sel_p(E/\Q)$ contains $E(\Q)\otimes\Qp/\Zp\cong(\Qp/\Zp)^r$ with $r\ge2$, so by Theorem~\ref{thm:PR}(iii) the $\Zp$-module $X^+/TX^+$ has rank at least $2$. Let $\fP=(T)$, a height-one prime of $\Lambda$; the localisation $\Lambda_\fP$ is a discrete valuation ring with uniformiser $T$ and residue field $\Qp$. By the structure theorem, $X^+_\fP\cong\bigoplus_j\Lambda_\fP/(T^{e_j})$ and $\ord_T\charid(X^+)=\sum_je_j$. The number of summands equals $\dim_{\Qp}X^+_\fP/TX^+_\fP=\dim_{\Qp}(X^+/TX^+)\otimes_{\Zp}\Qp\ge2$, so $\sum_j e_j\ge2$. Now apply Theorem~\ref{thm:PR}(i).
\end{proof}

\begin{remark}
The lemma is an exact divisibility in $\Lambda$, not a congruence modulo $p$. It uses no analytic input; combined with \cite[(12)]{PollackRubin}, $L_p^+(0)=(p-1)L(E,1)/\Omega_E$, it gives back $L(E,1)=0$. The same argument applies to $L_p^-$.
\end{remark}

\subsection{The criterion}

\begin{proposition}\label{prop:criterion}
If $\Aplus=[T^2]L_p^+$ is a $p$-adic unit, then $X^+\cong\Zp^2$ with $T$ acting as $0$, $\Sel_p(E/\Q)=E(\Q)\otimes\Qp/\Zp\cong(\Qp/\Zp)^2$, $\operatorname{rank}E(\Q)=2$, and $\Sha(E/\Q)[p^\infty]=0$.
\end{proposition}

\begin{proof}
By Lemma~\ref{lem:order}, $L_p^+=T^2u$ with $u\in\Lambda$ and $u(0)=\Aplus\in\Zp^\times$; hence $u\in\Lambda^\times$ and $\charid(X^+)=(T^2)$ by Theorem~\ref{thm:PR}(i). Localising at $\fP=(T)$, the module $X^+_\fP$ has length $2$ and, by the proof of Lemma~\ref{lem:order}, at least two cyclic summands; so $X^+_\fP\cong(\Lambda_\fP/T)^2$ and $TX^+_\fP=0$. At every height-one prime $\mathfrak{Q}\neq(T)$ the characteristic ideal is a unit, so $X^+_{\mathfrak Q}=0$. Hence the submodule $TX^+$ vanishes at every height-one prime, i.e.\ it is pseudo-null; since $\Lambda$ is a two-dimensional regular local ring, a finitely generated pseudo-null $\Lambda$-module is finite. By Theorem~\ref{thm:PR}(ii), $TX^+=0$. Thus $X^+$ is a finitely generated $\Zp$-module on which $T$ acts as $0$; its $\Zp$-torsion is a finite $\Lambda$-submodule, hence zero, so $X^+\cong\Zp^s$, and $\charid_\Lambda(\Zp^s)=(T^s)$ forces $s=2$. By Theorem~\ref{thm:PR}(iii), $\Hom(\Sel_p(E/\Q),\Qp/\Zp)\cong X^+/TX^+=X^+\cong\Zp^2$, i.e.\ $\Sel_p(E/\Q)\cong(\Qp/\Zp)^2$. It contains the divisible subgroup $E(\Q)\otimes\Qp/\Zp$ of corank $r\ge2$, so $r=2$ and the quotient $\Sha(E/\Q)[p^\infty]$ is divisible of corank $0$, hence finite and divisible, hence $0$.
\end{proof}

\begin{remark}\label{rem:criterion}
(a) Conversely, given Theorem~\ref{thm:PR} and $\operatorname{rank}E(\Q)\ge2$, the condition $X^+\cong\Zp^2$ (with $T=0$) implies $\charid(X^+)=(T^2)$ and hence $\Aplus\in\Zp^\times$. So the hypothesis of Proposition~\ref{prop:criterion} is equivalent to $X^+\cong\Zp^2$, which is strictly stronger than $\Sha(E/\Q)[p^\infty]=0$: for example $\Lambda/(T)\oplus\Lambda/(Tg)$, with $g$ a non-unit prime to $T$, has quotient by $T$ isomorphic to $\Zp^2$ but characteristic ideal different from $(T^2)$. If $\Aplus\equiv0\pmod p$, the criterion says nothing; this does not happen for $p<30000$.

(b) The use of Theorem~\ref{thm:PR}(ii) cannot be omitted: the module $(\Lambda/T)^2\oplus\Lambda/(p,T)$ has characteristic ideal $(T^2)$ and rank-two quotient by $T$, but would correspond to $\Sha(E/\Q)[p]\ne0$.

(c) Only the plus sign is used.
\end{remark}

\subsection{The coefficient \texorpdfstring{$\Aplus$}{A+} from the first layer}

\begin{lemma}\label{lem:finite-layer}
Put
\[
\Theta_1=\sum_{a\in(\Z/p^2\Z)^\times}\msp{a/p^2}\,(1+T)^{j(a)}\in\Q[T],
\]
with $j(a)$ represented in $\{0,\dots,p-1\}$. If all $\msp{a/p^2}$ lie in $\Zp$, then
\[
L_p^+\equiv-\Theta_1\pmod{\Phi_p(1+T)\Lambda},
\]
where $\Phi_p$ is the $p$-th cyclotomic polynomial, and consequently
\begin{equation}\label{eq:Aplus}
\Aplus\equiv-\sum_{a\in(\Z/p^2\Z)^\times}\msp{a/p^2}\binom{j(a)}{2}\pmod p .
\end{equation}
\end{lemma}

\begin{proof}
Let $\chi$ be a character of $\Gamma$ of order $p$. Viewed as a Dirichlet character it is even and primitive of conductor $p^2$, and $\chi(a)=\zeta^{j(a)}$ with $\zeta=\chi(\gamma)$. By \cite[(10)]{PollackRubin} with $n=1$ (here $\omega_1^+=1$ is an empty product and the sign is $(-1)^{(1+1)/2}=-1$),
\[
\chi(L_p^+)=-\frac{\tau(\chi)L(E,\bar\chi,1)}{\Omega_E},\qquad \tau(\chi)=\sum_{a\bmod p^2}\chi(a)e^{2\pi ia/p^2}.
\]
Here $\tau(\chi)$ is the Gauss sum of the primitive character $\chi$ of conductor $p^2$.\footnote{In \cite[p.~453]{PollackRubin} the Gauss sum of a character of order $p^k$ of $\Gamma$ is written with the root of unity $\zeta_{p^k}$; such a character has conductor $p^{k+1}$, and we use the Gauss sum of that conductor, as displayed. A different choice of the compatible system of roots of unity multiplies $L_p^\pm$ by an element of $\Gamma$, i.e.\ by $(1+T)^c$ for some $c\in\Zp$, which does not change $[T^2]L_p^+$ because $T^2\mid L_p^+$ (Lemma~\ref{lem:order}).}
Birch's formula: since $\sum_a\chi(a)e^{2\pi ina/p^2}=\bar\chi(n)\tau(\chi)$ for primitive $\chi$, we have $\sum_a\chi(a)f(z+a/p^2)=\tau(\chi)f_{\bar\chi}(z)$, where $f_{\bar\chi}=\sum\bar\chi(n)a_nq^n$; integrating from $i\infty$ to $0$ gives
\[
\tau(\chi)L(E,\bar\chi,1)=\sum_a\chi(a)\cdot2\pi i\int_{i\infty}^{a/p^2}f(z)\,dz .
\]
As $\chi$ is even, $\chi(-a)=\chi(a)$ and $[-r]^\pm=\pm[r]^\pm$, so the minus parts cancel and $\tau(\chi)L(E,\bar\chi,1)/\Omega_E=\sum_a\chi(a)\msp{a/p^2}=\chi(\Theta_1)$. Hence $L_p^++\Theta_1\in\Lambda$ vanishes at all $p-1$ roots of the distinguished polynomial $\Phi_p(1+T)$, which therefore divides it in $\Lambda$ by the Weierstrass division theorem. Modulo $p$ we have $\Phi_p(1+T)\equiv T^{p-1}$, so the coefficients of $T^k$, $k\le p-2$, of $L_p^+$ and $-\Theta_1$ agree modulo $p$; for $p\ge5$ this includes $k=2$, and $[T^2](1+T)^j=\binom j2$.
\end{proof}

\begin{remark}\label{rem:normalisation}
(a) The lemma uses only characters of exact order $p$; it does not use $L(E,1)=0$ or the value at the trivial character. As a consistency check, the same computation gives $[T^0]$ and $[T^1]$ of $L_p^+$ modulo $p$, which must vanish by Lemma~\ref{lem:order}; they do at every prime computed.

(b) In the computation we use $k(a)=\bigl((a^{p-1}-1)/p\bigr)\bmod p$, so that $a^{p-1}\equiv(1+p)^{k(a)}\pmod{p^2}$; since $\langle a\rangle^{p-1}=a^{p-1}$ we have $k(a)\equiv-j(a)$. Because $\binom{-k}{2}=\binom k2+k$ and $\sum_a\msp{a/p^2}k(a)\equiv0\pmod p$ (this is the vanishing of the linear term), \eqref{eq:Aplus} holds equally with $k(a)$ in place of $j(a)$.

(c) Multiplying $\Omega_E$ by a $p$-adic unit multiplies $\Aplus$ by a unit, and so does a change of $\gamma$; the unit property is insensitive to both.
\end{remark}

\subsection{Proof of Theorem A, given the computation}

\begin{corollary}\label{cor:theoremA}
Let $p\equiv3\pmod4$, $7\le p<30000$. Assume that the value of the right-hand side of \eqref{eq:Aplus} is nonzero modulo $p$ (this is established by the computation of Section~\ref{sec:computation}). Then $X^+\cong\Zp^2$, $\Sha(E/\Q)[p^\infty]=0$ and $\Sha(E/\Q(i))[p^\infty]=0$.
\end{corollary}

\begin{proof}
All values $\msp r$, $r\in\Q$, lie in $\frac12\Z$ (Section~\ref{subsec:exact}), so Lemma~\ref{lem:finite-layer} applies and gives $\Aplus\in\Zp^\times$; Proposition~\ref{prop:criterion} gives the statements over $\Q$. For $K=\Q(i)$: since $p$ is odd and $\Gal(K/\Q)$ has order $2$, restriction identifies $H^1(\Q,E[p^\infty])$ and $H^1(\Q,E^{(-1)}[p^\infty])$ with the two eigenspaces of $\Gal(K/\Q)$ acting on $H^1(K,E[p^\infty])$, compatibly with the local conditions and the Kummer maps; here $E^{(-1)}$ is the quadratic twist by $K/\Q$. Hence $\Sha(E/K)[p^\infty]\cong\Sha(E/\Q)[p^\infty]\oplus\Sha(E^{(-1)}/\Q)[p^\infty]$, and $E^{(-1)}\cong E$ over $\Q$ (Section~\ref{sec:curve}).
\end{proof}

\begin{remark}
In the proof of Theorem~\ref{thm:PR}(iii), the local analysis at the bad primes $v\in\{2,3,13\}$ requires $E(\Q_{\infty,w})[p]=0$ for $w\mid v$. This holds for $p\ge5$: $E$ has additive reduction of type III at $v$, the component group has order $2$ (over an unramified extension, at most $4$), the group of non-singular points of the reduction is an extension of the additive group of the residue field (a $v$-group) by a pro-$v$ formal group, and $\Q_{\infty,w}/\Q_v$ is unramified. Tamagawa numbers do not appear in the argument.
\end{remark}

\section{The computation}\label{sec:computation}

\subsection{Exact modular symbols}\label{subsec:exact}
The plus modular symbols were taken from eclib \cite{eclib} through SageMath, as
\begin{center}
\texttt{EllipticCurve([0,0,0,39,0]).modular\_symbol(sign=+1, implementation='eclib')}.
\end{center} eclib computes the plus eigenspace of modular symbols of level $N$ by exact linear algebra over $\Q$; the resulting map $r\mapsto\msp r$ is therefore $\lambda\cdot\varphi(r)$ with $\varphi$ exactly rational and a single global constant $\lambda\in\Q^\times$, which eclib fixes by a numerical comparison of periods. We checked $\lambda=1$ independently of eclib (Section~\ref{subsec:checks}).

From this map we built a table of the values $2(\msp{b/d}-\msp{a/c})$ for one matrix $\begin{psmallmatrix}a&b\\c&d\end{psmallmatrix}\in SL_2(\Z)$ over each of the $104832=[SL_2(\Z):\Gamma_0(48672)]$ points $(c:d)$ of $\mathbf P^1(\Z/48672\Z)$ (enumerated as $\mathbf P^1(\Z/32)\times\mathbf P^1(\Z/9)\times\mathbf P^1(\Z/169)$, of sizes $48\cdot12\cdot182$). The program checks that $\msp0=0$, that every class is filled exactly once, that all entries are integers, and that they lie in $[-8,8]$; building the table took $4$ seconds. By Manin's continued-fraction algorithm every $\msp{a/m}$ is a sum of table entries divided by $2$, so all plus modular symbols lie in $\frac12\Z$, and for odd $p$ they are $p$-integral. Hence \eqref{eq:Aplus} is evaluated by exact integer arithmetic, without rounding.

\subsection{The sum}
For $q=p^2$ and $a/q$ with continued-fraction denominators $q_{-1}=0,q_0=1,\dots,q_n=q$, one has $2\msp{a/q}=\sum_i e(q_{i-1},q_i)$, where $e(\cdot,\cdot)$ is the table entry of the corresponding class; the sign of the determinant is irrelevant because $\msp{-r}=\msp r$. Each symbol costs $O(\log q)$ table look-ups. We accumulate
\[
s_0=\sum 2\msp{a/q},\quad s_1=\sum2\msp{a/q}k(a),\quad s_2=\sum2\msp{a/q}\binom{k(a)}2
\]
over $1\le a\le q/2$, $p\nmid a$ (the terms for $a$ and $q-a$ coincide, which absorbs the factor $2$), reducing modulo $p$ periodically. The program asserts $s_0=0$ exactly and $s_1\equiv0\pmod p$, and outputs $\Aplus\equiv-s_2\pmod p$. Three implementations of the summation were used: a direct one for $7\le p<3000$ ($217$ primes); a ``folded'' one for $3000\le p<10000$ ($401$ primes), which reads the edges above $\sqrt q$ from the Euclidean remainder sequence of the inverse residue so that each edge is counted once; and the folded one with cached common continued-fraction prefixes for $10000\le p<30000$ ($1014$ primes). The cost is $O(p^2\log p)$ integer operations per prime. The total running time was $4019$ seconds on a single core of an arm64 machine running macOS (largest single prime: $p=29983$, $7.0$ seconds), with peak memory below $0.5$ GB.

\subsection{Result}
For every one of the $1632$ inert primes $7\le p\le29983$ the computed residue $\Aplus\bmod p$ lies in $\{1,\dots,p-1\}$; it is never $0$. The first values are
\[
\begin{array}{c|ccccccc}
p&7&11&19&23&31&43&47\\\hline
\Aplus\bmod p&6&8&9&15&6&27&32
\end{array}
\]
and the complete list is in the ancillary file. By Corollary~\ref{cor:theoremA} this proves Theorem~A, subject to the correctness of the computation.

\subsection{Checks}\label{subsec:checks}
All recomputations listed below agree with the main computation; no discrepancy was found.

\begin{enumerate}
\item \emph{Different summation paths.} $102$ primes were computed by at least two of the three summation modes (or by one mode with different parameters).
\item \emph{Group-ring reconstruction.} For $20$ primes between $751$ and $29123$, chosen by a pseudo-random rule fixed in advance, the coefficients of $\Theta_1$ were rebuilt directly in the group ring from Teichm\"uller representatives and powers of $1+p$, and $\Aplus$ was recomputed.
\item \emph{Independent rebuilding from the eclib object.} Two further programs, which do not use our table or our continued-fraction code but do use the same eclib modular-symbol object, recomputed $\Aplus$ at all $43$ inert primes $p\le443$ and at $24$ further primes between $503$ and $29303$ ($67$ primes in all; together with item~2, $86$ distinct primes have a recomputation of $\Aplus$ by a second program).
\item \emph{Checks independent of eclib.} (a) At all $14$ inert primes $7\le p\le107$, $53528$ plus modular symbols of the relevant layers were computed by certified (interval) numerical integration and agree with ours one by one. (b) At the $12$ inert primes $7\le p\le83$ we computed the twisted central values $L(E,\bar\chi,1)$ for all characters $\chi$ of order $p$ of $(\Z/p^2\Z)^\times$ with PARI's \texttt{lfuntwist}, recovered the integers $\sum_{k(a)=k}2\msp{a/p^2}$ by an inverse discrete Fourier transform (rounding error below $10^{-33}$), and found that these integer vectors coincide componentwise with those obtained from our table. (c) At $p=7,11,19,23,31,43,47$ a third implementation based on twisted $L$-values with independently computed Fourier coefficients reproduced $\Aplus$. (d) At $20$ further inert primes between $2179$ and $27191$, drawn without replacement from all $1632$ by a pseudo-random rule with a seed fixed in advance, $\Aplus\bmod p$ was recomputed without the eclib object: the plus modular symbols of level $N$ were rebuilt from the Fourier coefficients of the newform and the real period of the minimal model, by evaluating $8520$ cusp integrals numerically with explicit error bounds, rounding them to half-integers, and solving the Manin relations exactly over $\Q$; all $20$ residues agree with ours. These checks also pin down the global constant $\lambda=1$: a different constant would rescale all the integer vectors in (b).
\item \emph{Another code path.} At $p=7,11,19,23$, SageMath's \texttt{padic\_lseries} (supersingular case, using its own Riemann sums) gives $[T^2]L_\alpha=(\Aplus+\alpha A_-)/p$ with residues agreeing with ours (cf.~\eqref{eq:Lalpha} below).
\end{enumerate}

\begin{remark}[what the checks do not cover]\label{rem:shared}
All $1632$ values of $\Aplus$ rest on a single modular-symbol input, the eclib object of level $48672$. For $p\le107$ there are checks that do not use eclib (item~4(a)--(c)); for $p>107$ the recomputations in items~1--3 use different code but the same eclib object, and the only recomputation that does not use it is the random sample of $20$ primes in item~4(d). The only check over the whole range that does not use modular symbols at all is the agreement with $-2d_p$ in Observation~C; since that agreement is predicted by a conjecture, it is a consistency check and not a second proof. The step from ``$\lambda=1$ at small primes'' to ``$\lambda=1$ everywhere'' relies on the fact that eclib uses one rational scaling for all symbols of the given level.
\end{remark}

\section{The primes \texorpdfstring{$2$, $3$ and $13$}{2, 3 and 13}}\label{sec:bad}

The arguments in this section are descents of a standard kind; we include them for completeness. None depends on the generalised Riemann hypothesis.

\subsection{\texorpdfstring{$p=2$}{p=2}}

\begin{proposition}\label{prop:p2}
$\Sha(E/\Q)[2^\infty]=0$.
\end{proposition}

\begin{proof}
Let $\varphi\colon E\to E'\colon y^2=x^3-156x$ be the $2$-isogeny with kernel $\{O,(0,0)\}$ and $\hat\varphi$ its dual. A complete $2$-isogeny descent (see \cite[\S X.4]{Silverman} and \cite[\S3.6]{Cremona97}), carried out with mwrank \cite{eclib}, shows that the Selmer groups $\Sel^{(\varphi)}(E/\Q)$ and $\Sel^{(\hat\varphi)}(E'/\Q)$ both have $\F_2$-dimension $2$ and are spanned by the images of $E'(\Q)$ and $E(\Q)$ respectively. Hence $\Sha(E/\Q)[\varphi]=0$ and $\Sha(E'/\Q)[\hat\varphi]=0$. Since $\hat\varphi\circ\varphi=[2]$, $\varphi$ maps $\Sha(E/\Q)[2]$ into $\Sha(E'/\Q)[\hat\varphi]=0$, so $\Sha(E/\Q)[2]\subseteq\Sha(E/\Q)[\varphi]=0$, and therefore $\Sha(E/\Q)[2^\infty]=0$.
\end{proof}

mwrank reports $\#\Sha(E/\Q)[2]=1$ and a $2$-Selmer group of dimension $3=2+\dim E(\Q)[2]$, with the result flagged as certain. This is a finite computation (local solubility tests and a search for rational points) and uses neither the generalised Riemann hypothesis nor BSD.

\subsection{\texorpdfstring{$p=3$}{p=3}}

Let $\theta$ be a root of $X^8+78X^4-507$ and $L=\Q(i,\theta)$. The following facts were computed with PARI/GP 2.17.3: $x^4+78x^2-507$ (a factor of the $3$-division polynomial $3(x^4+78x^2-507)$) is irreducible; $[L:\Q]=16$ and $L/\Q$ is Galois with group the semidihedral group $SD_{16}$, $\Gal(L/K)\cong C_8$; all eight nonzero $3$-torsion points of $E$ are defined over $L$ and $L=K(E[3])$; $i\notin F_8:=\Q(\theta)$; $\Q(\zeta_3)\subset L$ but $\Q(\zeta_9)\not\subset L$; $3$ has a single place in $L$; and the class group of $F_8$ is $Cl(F_8)\cong\Z/12\Z\times\Z/2\Z$, certified unconditionally by \texttt{bnfcertify}.

\begin{proposition}\label{prop:p3}
$\Sha(E/K)[3^\infty]=0$ and $\Sha(E/\Q)[3^\infty]=0$.
\end{proposition}

\begin{proof}
\emph{Set-up.} Since $E$ has CM by $\cO_K=\Z[i]$ defined over $K$ and $3$ is inert in $K$, $E[3]$ is free of rank one over $\cO_K/3=\F_9$ and $G_K$ acts through a character $\rho\colon G_K\to\F_9^\times$. By the computation above, $[L:K]=8$, so $\rho$ induces $\Gal(L/K)\cong\F_9^\times\cong C_8$, and $\Gal(L/\Q)\cong SD_{16}=\langle\sigma,\tau\mid\sigma^8=\tau^2=1,\ \tau\sigma\tau^{-1}=\sigma^3\rangle$ with $\langle\sigma\rangle=\Gal(L/K)$. In particular $E(K)[3]=0$.

\emph{Step 1 (lower bound).} $\operatorname{rank}E(K)=4$ and $E(K)[3]=0$, so $\dim_{\F_3}E(K)/3E(K)=4$ and $\dim\Sel_3(E/K)\ge4$.

\emph{Step 2 (local conditions away from $3$).} Let $v$ be a place of $K$ above $2$ or $13$. The reduction at $v$ is additive with Tamagawa number $c_v\in\{2,4\}$; the filtration $E(K_v)\supset E_0(K_v)\supset E_1(K_v)$ has graded pieces of order $c_v$, a power of the residue characteristic $\ell\in\{2,13\}$, and a pro-$\ell$ group. Hence $E(K_v)/3E(K_v)=0$ and every Selmer class is locally trivial at $v$. At places not dividing $78$, Selmer classes are unramified.

\emph{Step 3 (passage to $L$).} Since $[L:K]=8$ is prime to $3$, restriction gives an injection $\Sel_3(E/K)\hookrightarrow\Hom(G_L,E[3])^{\Gal(L/K)}$. As $\mu_3\subset L$ and $E[3]$ is a trivial $G_L$-module, $\Hom(G_L,E[3])\cong(L^\times/L^{\times3})\otimes_{\F_3}M$ with $M=\Hom(\mu_3,E[3])$, equivariantly for $\Gal(L/K)$. By Step~2 the image consists of classes $x$ such that $L(x^{1/3})/L$ is unramified outside the place above $3$; it lies in $(L(S,3)\otimes M)^{C_8}$, where $S$ is the set of places of $L$ above $3$ and $L(S,3)=\{x\in L^\times/L^{\times3}: v_w(x)\equiv0\ (3)\text{ for all }w\notin S\}$. The representation $M$ is the twist of the faithful $\F_9$-character $\rho$ by the quadratic character of $\mu_3$, hence faithful; $M^*=\Hom(M,\F_3)$ is an irreducible $\F_3[C_8]$-module with $\operatorname{End}_{C_8}(M^*)\cong\F_9$. As $\F_3[C_8]$ is semisimple, $\dim_{\F_3}(L(S,3)\otimes M)^{C_8}=2m$, where $m$ is the multiplicity of $M^*$ in $L(S,3)$.

\emph{Step 4 (units).} There is an exact sequence $0\to\cO_{L,S}^\times/3\to L(S,3)\to Cl_S(L)[3]\to0$. Since $K$ has a single (complex) archimedean place and $S$ consists of a single place fixed by $\Gal(L/K)$, Herbrand's theorem gives $\cO_{L,S}^\times\otimes\Q\cong\Q[C_8]$; as $3\nmid8$, the torsion-free quotient modulo $3$ is $\cong\F_3[C_8]$, which contains $M^*$ exactly once. The torsion contributes $\mu(L)/3=\mu_3$ (as $\mu_9\not\subset L$), on which $C_8$ acts through a character of order $2$, not by $M^*$. So $M^*$ occurs exactly once in $\cO_{L,S}^\times/3$.

\emph{Step 5 (class group).} The class of the unique prime above $3$ is fixed by $\Gal(L/K)$, so the $M^*$-isotypic parts of $Cl_S(L)[3]$ and $Cl(L)[3]$ coincide. The two faithful irreducible $\F_3$-representations $V_1,V_2$ of $SD_{16}$ are $2$-dimensional and absolutely irreducible (their characters take values in $\Z[\sqrt{-2}]$, and $\sqrt{-2}\in\F_3$), and they restrict to $C_8$ as the two faithful irreducible $\F_3[C_8]$-modules; every non-faithful irreducible representation of $SD_{16}$ is trivial on $\sigma^4$ and contains no faithful constituent of $C_8$. Hence the multiplicity of $M^*$ in $Cl(L)[3]$ is at most $m_1+m_2$, where $m_j$ is the multiplicity of $V_j$ in $Cl(L)[3]$. Let $H=\Gal(L/F_8)=\langle c\rangle$. Since $i\notin F_8$, $c\notin\langle\sigma\rangle$, so $c$ is a non-central involution; its trace on $V_j$ is $0$, so $\dim V_j^{c}=1$. Since $[L:F_8]=2$ is prime to $3$, $Cl(F_8)[3]\cong Cl(L)[3]^H$. Therefore
\[
m_1+m_2\le\sum_V m_V\dim V^c=\dim_{\F_3}Cl(F_8)[3]=1 .
\]

\emph{Step 6 (conclusion).} Steps 3--5 give $m\le1+1=2$, hence $\dim\Sel_3(E/K)\le4$. With Step~1, $\Sel_3(E/K)=E(K)/3E(K)$, so $\Sha(E/K)[3]=0$ and $\Sha(E/K)[3^\infty]=0$. Finally $\mathrm{cor}\circ\mathrm{res}=[K:\Q]=2$ is invertible on $3$-primary groups, so restriction $\Sha(E/\Q)[3^\infty]\to\Sha(E/K)[3^\infty]$ is injective.
\end{proof}

The only computational input is PARI's certified class group of the degree-$8$ field $F_8$ (\texttt{bnfcertify} removes the dependence on the generalised Riemann hypothesis). The argument does not use the local conditions at $3$.

\subsection{\texorpdfstring{$p=13$}{p=13}}

Here $13=\pi\bar\pi$ with $\pi=3+2i$ splits in $K$, and $E$ has bad (additive) reduction at $13$. We use the following theorem of Rubin.

\begin{theorem}[{Rubin \cite[Thm.~3.3, p.~35]{Rubin91}}]\label{thm:rubin33}
Let $K$ be imaginary quadratic, $p$ a prime not dividing the number of roots of unity in the Hilbert class field $H$ of $K$, and $F\supseteq H$ a finite abelian extension of $K$ with $p\nmid[F:K]$. Let $A$ be the $p$-part of the ideal class group of $F$, $\cE$ the unit group of $F$, and $\cC$ the group of elliptic units of $F$ as defined in \cite[\S1, pp.~28--29]{Rubin91}. Then for every irreducible $\Zp$-representation $\chi$ of $\Gal(F/K)$, $\#A^\chi=\#\bigl(\Zp\otimes\cE/\cC\bigr)^\chi$.
\end{theorem}

For us $K=\Q(i)$, $H=K$, $\#\mu(K)=4$ is prime to $13$, and $F=K(E[\pi])$ is cyclic of degree $12$ over $K$; so the hypotheses hold, and all irreducible $\Z_{13}$-representations of $\Gal(F/K)$ are characters.

\begin{proposition}\label{prop:p13}
$\Sha(E/K)[13^\infty]=0$ and $\Sha(E/\Q)[13^\infty]=0$.
\end{proposition}

\begin{proof}
$E[\pi]$ is one-dimensional over $\cO_K/\pi=\F_{13}$; $G=\Gal(F/K)$ acts on it by a faithful character $\chi\colon G\cong C_{12}\to\F_{13}^\times$. Let $A=Cl(F)\otimes\Z_{13}$.

\emph{Step 1 (local conditions).} (i) For $v\mid2\cdot3$: $E(K_v)$ is an extension of a group of order $c_v\in\{2,4\}$ by $E_0(K_v)$, which is an extension of the additive group of the residue field by a pro-$\ell$ group ($\ell\in\{2,3\}$); so $E(K_v)[\pi]=0$ and, as $v\nmid13$, $E(K_v)/\pi E(K_v)=0$. (ii) At $\bar\pi$: $K_{\bar\pi}=\Q_{13}$ and $\pi$ is a unit there, so $[\pi]$ is an automorphism of the formal group and acts on $E_0/E_1\cong\F_{13}$ by the nonzero scalar $\pi\bmod\bar\pi$; the component group has order $2$; so $E(K_{\bar\pi})/\pi=0$. (iii) At $\pi$: $K_\pi=\Q_{13}$, $E$ has type III with $c=2$, $E(K_\pi)[\pi]=0$ (the decomposition group in $G$ of the place of $F$ above $\pi$ is all of $G$, see Step~2, and $\chi$ is nontrivial on it, so $G_{K_\pi}$ fixes no nonzero point of $E[\pi]$), and $[\pi]$ kills $E_0/E_1\cong\F_{13}$ while mapping $E_1$ onto $E_2$; hence $E(K_\pi)/\pi\cong E_0/E_1$ is one-dimensional, and $P=(3,12)$ maps to a nonzero element (its reduction modulo $13$ is a non-singular point other than $O$).

\emph{Step 2 (strict Selmer group).} Let $\Sel_0\subseteq\Sel_\pi(E/K)$ be the subgroup of classes that are also locally trivial at $\pi$. By Step 1(iii), $\Sel_\pi(E/K)=\Sel_0+\langle\kappa(P)\rangle$ and $\dim\Sel_\pi=\dim\Sel_0+1$. Since $13\nmid[F:K]$, restriction identifies $H^1(K,E[\pi])$ with $\Hom_G(G_F,E[\pi])$. By Step~1 and good reduction elsewhere, a class lies in $\Sel_0$ if and only if the corresponding homomorphism is unramified everywhere and trivial on the decomposition groups of the places of $F$ above $2$, $3$ and $13$. For such a place $w$ the decomposition group $D_w\subseteq G$ is nontrivial (the places of $K$ above $2$ and $3$ ramify in $F/K$, and $\pi$ and $\bar\pi$ each have a single place above them in $F$, so that $D_w=G$ there), and for $\sigma\in D_w$ a $G$-equivariant homomorphism $f$ on the class group satisfies $f([w])=f(\sigma[w])=\chi(\sigma)f([w])$; as $\chi$ is faithful, $f([w])=0$ automatically. Hence $\Sel_0\cong\Hom_G(A/13A,E[\pi])$, whose dimension is $\dim_{\F_{13}}(A/13A)^{\chi'}$ for one character $\chi'\in\{\chi,\chi^{-1}\}$ (which one depends on conventions; the argument does not need to know).

\emph{Step 3 (lower bound).} As an $\cO_K$-module $E(K)\cong\cO_K^2\oplus\Z/2$, and $E(K)[\pi]=0$, so $\dim E(K)/\pi E(K)=2\le\dim\Sel_\pi$. Thus $\dim(A/13A)^{\chi'}\ge1$, i.e.\ $\#A^{\chi'}\ge13$.

\emph{Step 4 (upper bound).} By Theorem~\ref{thm:rubin33}, $\#A^\psi=\#(\Z_{13}\otimes\cE/\cC)^\psi$ for $\psi\in\{\chi,\chi^{-1}\}$ (Teichm\"uller lifts). Since $\mu_{13}\not\subset F$, $\Z_{13}\otimes\cE$ is torsion free, and by Herbrand's theorem ($F$ is totally imaginary and $K$ has one complex place) each nontrivial $\psi$-component $(\Z_{13}\otimes\cE)^\psi$ is free of rank one over $\Z_{13}$. Let $\ff=(1+i)^3\cdot3\cdot13$ be the conductor of the Hecke character of $E$ (so that $4\,\mathrm N\ff=N$), $\tau$ a primitive $\ff$-division point, and $u=N_{K(\ff)/F}\Theta(\tau;L,(5))$, where $\Theta(z;L,(5))$ differs from Rubin's $\Theta_0(z;(5))^{12}$ \cite[p.~28]{Rubin91} by a constant in $K^\times$ and is, up to a constant in $K^\times$, $\Delta^{24}\psi_5(x)^{-12}$ with $\psi_5$ the $5$-division polynomial (cf.~\cite{deShalit}). Since $\pi\mid\ff$ and $\psi_E((\alpha))=\alpha\equiv1\pmod{\pi}$ whenever $\alpha\equiv1\pmod{\ff}$, where $\psi_E$ is the Hecke character of $E$, the field $F$ is contained in the ray class field $K(\ff)$, so the norm $N_{K(\ff)/F}$ makes sense and $FK(\ff)=K(\ff)$. The map $\cO_K^\times\to(\cO_K/\ff)^\times$ is injective and $(5)$ is prime to $6\ff$, so by definition \cite[\S1, pp.~28--29]{Rubin91} the group $\cC$ contains $\alpha_{\tau,(5)}(1)^{\sigma-1}$ for every $\sigma\in G$, where $\alpha_{\tau,(5)}(1)=N_{K(\ff)/F}\,\Theta_0(\tau;(5))$; and $u=c\cdot\alpha_{\tau,(5)}(1)^{12}$ with $c\in K^\times$. For $\psi\ne1$ the idempotent $e_\psi$ kills $c$, the factor $12$ is a unit in $\Z_{13}$, and $\psi(\sigma)-1\in\Z_{13}^\times$ for a generator $\sigma$ of $G$; hence
\[
\Z_{13}\,e_\psi u\subseteq e_\psi(\Z_{13}\otimes\cC)\subseteq(\Z_{13}\otimes\cE)^\psi,\qquad\text{so}\qquad\#\bigl(\Z_{13}\otimes\cE/\cC\bigr)^\psi\le\bigl[(\Z_{13}\otimes\cE)^\psi:\Z_{13}\,e_\psi u\bigr].
\]
(We do not claim that $e_\psi u$ generates $e_\psi(\Z_{13}\otimes\cC)$; $\cC$ also contains elliptic units of other conductors, and only the inclusion is used.) For the three primes $\fq$ of $K$ of degree one with $\mathrm N\fq=q\in\{24337,608401,876097\}$, which split completely in $K(\ff)$ and satisfy $q\equiv1\pmod{169}$, we computed the reduction of $e_\psi u$ modulo a prime above $\fq$ (exact arithmetic in $\F_q$, with the idempotent taken modulo $169$) and the order of its $((q-1)/169)$-th power: it is $13$ for one of the two characters and $169$ for the other, at all three $q$. Since $v_{13}(q-1)=2$ for each of the three $q$, the image of $e_\psi u$ in $\F_q^\times\otimes\Z_{13}\cong\Z/169$ thus has order $169$ for one character and order $13$ for the other; a $13$th (resp.\ $169$th) power in $(\Z_{13}\otimes\cE)^\psi$ would map to a $13$th (resp.\ $169$th) power. As $(\Z_{13}\otimes\cE)^\psi\cong\Z_{13}$, the index $[(\Z_{13}\otimes\cE)^\psi:\Z_{13}e_\psi u]$ is therefore $1$ for one character and at most $13$ for the other, and by the displayed inequality the same bounds hold for $\#(\Z_{13}\otimes\cE/\cC)^\psi$. By Step~3 the $\chi'$-component of $A$ is nontrivial, so $\chi'$ is the second character and $\#A^{\chi'}=13$, i.e.\ $A^{\chi'}\cong\Z/13$.

\emph{Step 5 (conclusion).} $\dim\Sel_\pi(E/K)=1+1=2=\dim E(K)/\pi E(K)$, so $\Sha(E/K)[\pi]=0$ and hence $\Sha(E/K)[\pi^\infty]=0$. Complex conjugation interchanges the $\pi$- and $\bar\pi$-primary parts ($E$ is defined over $\Q$), so $\Sha(E/K)[13^\infty]=0$; and restriction $\Sha(E/\Q)[13^\infty]\hookrightarrow\Sha(E/K)[13^\infty]$ is injective since $\mathrm{cor}\circ\mathrm{res}=2$.
\end{proof}

The power-residue computation in Step~4 is exact arithmetic in finite fields; no class group or unit group of a number field is computed, and the generalised Riemann hypothesis is not used.

\section{The ordinary primes}\label{sec:split}

Coates, Liang and Sujatha proved the following for five curves $E_i\colon y^2=x^3-D_ix$ with $D\in\{-14,17,-33,-34,-39\}$; the fifth, $D_5=-39$, is our curve.

\begin{theorem}[{Coates--Liang--Sujatha \cite[Thm.~1.3]{CLS2}}]\label{thm:CLS}
For each of the five curves and each prime $p\equiv1\pmod4$ of good reduction with $p<30000$, $\Sha(E_i/\Q)[p^\infty]$ is finite, and it is zero except possibly in the four cases $(E_1,29)$, $(E_1,277)$, $(E_4,577)$, $(E_5,17)$.
\end{theorem}

The proof applies the criterion of \cite[Thm.~2.2]{CLS1} (which rests on Rubin's main conjecture \cite{Rubin91}) to a quantity $c_p^+$ computed exactly by the authors in MAGMA. For our curve the hypotheses of \cite[Thm.~2.2]{CLS1} hold at every $p\equiv1\pmod4$, $p\ne13$: $p$ is odd and prime to the number $w=4$ of roots of unity of $K$ and to $\alpha(E)=1+i$; $p$ splits in $K$; the root number is $+1$, so the analytic and algebraic ranks have the same parity; and $\#\tilde E(\F_p)=p+1-a_p$ is prime to $p$ (it is even and $|a_p|\le2\sqrt p$). Only the ranges $p<1000$ and $29000<p<30000$ of the computation are tabulated in \cite{CLS2}; we have not repeated it. At $p=17$ the quantity $c_{17}^+$ has $17$-adic valuation $3$ (its factorisation is given at the end of \cite{CLS2}, p.~23 of the arXiv version), which by \cite[Thm.~2.2]{CLS1} gives finiteness only. The paragraph after \cite[Thm.~1.3]{CLS2} reports that $p$-adic height computations of C.~Wuthrich confirm the vanishing at $17$. A proof of the vanishing at $17$ is given in Banwait's preprint:

\begin{theorem}[{Banwait \cite[Prop.~7.4 and Cor.~7.5]{Banwait}, preprint}]\label{thm:Banwait}
$\Sha(E/\Q)[17^\infty]=0$. Consequently $\Sha(E/\Q)[p^\infty]=0$ for every prime $p\equiv1\pmod 4$ with $p\ne13$ and $p<30000$.
\end{theorem}

Banwait's argument applies his criterion \cite[Thm.~3.8]{Banwait} (rank two, CM by a maximal order, a split good prime $p\ge5$ outside an explicit finite set, and a unit condition on a second derivative) and verifies the unit condition at $17$ using eclib modular symbols with the error bounds of \cite[\S3]{SteinWuthrich}. (Banwait notes that \cite{SteinWuthrich} assumes throughout that the curve has no CM, and argues that the error bounds of \cite[\S3]{SteinWuthrich} do not use this assumption; we have not checked this independently.) We recomputed with SageMath that the coefficient of $T^2$ in the $17$-adic $L$-series of $E$ is $13+O(17)$, a $17$-adic unit, but SageMath also uses eclib modular symbols, so this is not independent of \cite{Banwait}.

\begin{remark}
Banwait's unit criterion fails for this curve at $p=5$ and at $p=15289$, where the relevant regulator is not a unit; at these primes the vanishing comes from Theorem~\ref{thm:CLS}. So the statement for the ordinary primes rests on \cite{CLS2} for all $p\ne17$, and on \cite{Banwait} at $p=17$.
\end{remark}

\section{Proof of Theorem B}\label{sec:combined}

There are $3245$ primes below $30000$, and every one of them falls into exactly one row of Table~\ref{tab:sources}: the bad primes $2,3,13$ (Propositions~\ref{prop:p2}, \ref{prop:p3}, \ref{prop:p13}); the $1610$ primes $p\equiv1\pmod 4$, $p\ne13$ (Theorems~\ref{thm:CLS} and~\ref{thm:Banwait}); and the $1632$ primes $p\equiv3\pmod 4$, $p\ge7$ (Theorem~A). Indeed $3245=3+1610+1632$. \qed

\smallskip
The strength of Theorem~B is that of its weakest input: the unpublished intermediate range of the computation behind \cite[Thm.~1.3]{CLS2}, the preprint \cite{Banwait} at $p=17$, and the computations of this paper.

\section{The leading-term congruence at supersingular primes: numerical evidence}\label{sec:leading}

This section and the next are not used in the proofs above.

\subsection{Set-up}
Let $E/\Q$ be given by a global minimal model with N\'eron differential $\omega=dx/2y$, and let $p\ge5$ be a prime of good supersingular reduction with $a_p=0$. Let $\Dp(E)=\Qp\otimes H^1_{\rm dR}(E/\Q)$, with basis $\omega,\eta=x\omega$, and let $\varphi$ be the Frobenius on $\Dp(E)$ normalised as in \cite[\S3.5]{SteinWuthrich}; its characteristic polynomial is $X^2-(a_p/p)X+1/p$, so $\varphi^2=-1/p$. We write $F=p\varphi$ (the crystalline Frobenius of Kedlaya's algorithm and of PARI).

Let $\alpha$ be a root of $X^2+p$, and let $L_\alpha(E,T)$ be the $p$-adic $L$-series of Amice--V\'elu and Vishik, normalised so that $L_\alpha(E,0)=(1-1/\alpha)^2L(E,1)/\Omega_E$. Pollack \cite[(1) and Lemma~4.1]{Pollack} showed
\begin{equation}\label{eq:Lalpha}
L_\alpha(E,T)=L_p^+(E,T)\log_p^+(T)+\alpha\,L_p^-(E,T)\log_p^-(T),\qquad \log_p^\pm(0)=1/p,
\end{equation}
with $L_p^\pm\in\Lambda$; the functions $L_p^\pm$ are those of Section~\ref{sec:curve} (\cite[\S7]{PollackRubin}; one checks that \eqref{eq:Lalpha} at $T=0$ is consistent with $\chi_0(L_p^+)=(p-1)L(E,1)/\Omega_E$ and $\chi_0(L_p^-)=2L(E,1)/\Omega_E$ from \cite[(12)]{PollackRubin}, using $1/\alpha=-\alpha/p$). Writing $L_\alpha=G+\alpha H$, the $\Dp(E)$-valued $p$-adic $L$-series of Perrin-Riou is \cite[\S3.5]{SteinWuthrich}
\[
\mathbf L_p(E,T)=G\,\omega+a_pH\,\omega-pH\,\varphi(\omega)=L_p^+\log_p^+\cdot\omega-p\,L_p^-\log_p^-\cdot\varphi(\omega)\qquad(a_p=0).
\]

For $r=\operatorname{rank}E(\Q)\ge1$, let $\Reg_p\in\Dp(E)$ be the $p$-adic regulator of \cite[Lemma~4.2]{SteinWuthrich} (a corrected form of \cite[Lemma~2.6]{PR03}), computed on a basis of $E(\Q)/E(\Q)_{\rm tors}$, and $\Reg_\gamma=\Reg_p/\log_p(\kappa(\gamma))^r$ \cite[\S4.4]{SteinWuthrich}; with our choice of $\gamma$, $\log_p\kappa(\gamma)=\ell:=\log_p(1+p)$. Write $\ell=p\ell'$, so $\ell'\equiv1\pmod p$.

\begin{definition}
Write $\Reg_p=R_1\,\omega+R_2\,F\omega$ with $R_1,R_2\in\Qp$. When $R_1,R_2\in p^2\Zp$, put $d_p=(R_1/p^2)\bmod p$ and $e_p=(R_2/p^2)\bmod p$.
\end{definition}

PARI's \texttt{ellpadic\-regulator} returns $\Reg_p$ in the basis $(\omega,F\omega)$; SageMath's \texttt{padic\_\-regulator} returns it in the basis $(\omega,\varphi\omega)$, and the two agree after the change of basis at all primes where we compared them. In every case computed below, $R_1,R_2\in p^2\Zp$; we have not written a general proof of this.

\subsection{The conjecture}

\begin{conjecture}[{Bernardi--Perrin-Riou; in the form of \cite[Conj.~5.2, p.~1775]{SteinWuthrich}; cf.\ \cite{BPR93,PR03,Sprung}}]\label{conj:BPR}
The order of vanishing of $\mathbf L_p(E,T)$ at $T=0$ equals $r=\operatorname{rank}E(\Q)$, and its leading coefficient $\mathbf L_p^*(E,0)$ satisfies
\begin{equation}\label{eq:BPR}
(1-\varphi)^{-2}\,\mathbf L_p^*(E,0)=\frac{\prod_\ell c_\ell\cdot\#\Sha(E/\Q)}{\#E(\Q)_{\rm tors}^2}\,\Reg_\gamma(E/\Q)\quad\text{in }\Dp(E).
\end{equation}
\end{conjecture}

We use the normalisation of \cite{SteinWuthrich}; the regulator in \cite[\S1.3]{PR03} is normalised differently for $r\ge2$ (see the footnote to \cite[Lemma~4.2]{SteinWuthrich}). We have not seen the original note \cite{BPR93}.

\subsection{From the conjecture to a congruence}
Put $C=\prod_\ell c_\ell/\#E(\Q)_{\rm tors}^2$.

\begin{proposition}\label{prop:derivation}
Suppose $r=2$, $a_p=0$, $T^2\mid L_p^\pm$, $R_1,R_2\in p^2\Zp$, $C\cdot\#\Sha(E/\Q)\in\Zp$ (with $\Sha(E/\Q)$ finite), and that \eqref{eq:BPR} holds. Write $a_\pm=[T^2]L_p^\pm$, $u=R_1/p^2$ and $e'=R_2/p^2$. Then
\[
a_+=C\cdot\#\Sha\cdot\frac{(p-1)u+2p\,e'}{\ell'^2},\qquad -a_-=C\cdot\#\Sha\cdot\frac{(p-1)e'-2u}{\ell'^2},
\]
and in particular
\[
\Aplus\equiv-C\cdot\#\Sha\cdot d_p\pmod p .
\]
\end{proposition}

\begin{proof}
Taking the coefficient of $T^2$ in $\mathbf L_p$ and using $\log_p^\pm(0)=1/p$ gives $\mathbf L_p^*(E,0)=(a_+/p)\,\omega-a_-\,\varphi\omega$. On the other side, $\Reg_p=R_1\,\omega+pR_2\,\varphi\omega$, since $F=p\varphi$. Since $\varphi^2=-1/p$ we have $(1-\varphi)^2=(1-1/p)-2\varphi$, and
\[
(1-\varphi)^2\Reg_p=\Bigl[(1-\tfrac1p)R_1+2R_2\Bigr]\omega+\Bigl[(p-1)R_2-2R_1\Bigr]\varphi\omega .
\]
Multiplying \eqref{eq:BPR} by $(1-\varphi)^2$ and comparing coordinates, with $\Reg_\gamma=\Reg_p/\ell^2$:
\[
a_+=C\#\Sha\,\frac{(p-1)R_1+2pR_2}{\ell^2},\qquad -a_-=C\#\Sha\,\frac{(p-1)R_2-2R_1}{\ell^2}.
\]
Substituting $R_1=p^2u$, $R_2=p^2e'$, $\ell=p\ell'$ gives the displayed formulas. Reducing modulo $p$, with $\ell'\equiv1$ and $p-1\equiv-1$, gives the congruence.
\end{proof}

The minus sign comes from the factor $(1-1/p)$ in the $\omega$-coordinate of $(1-\varphi)^2$. If, moreover, the complex BSD formula holds, then $C\cdot\#\Sha=q_\infty:=L''(E,1)/(2\Omega_E\Reg_\infty)$, and Proposition~\ref{prop:derivation} predicts
\begin{equation}\label{eq:prediction}
\Aplus\equiv-q_\infty\,d_p\pmod p,\qquad a_+\equiv q_\infty\frac{(p-1)u+2pe'}{\ell'^2}\pmod{p^2}.
\end{equation}
The second congruence requires $u\bmod p^2$, $e'\bmod p$ and $\ell'\bmod p^2$.

\begin{remark}
The test \eqref{eq:prediction} is sensitive to global normalisation errors (a factor $2$ in the period, Tamagawa numbers, the sign conventions), but not to the choice of a basis of $E(\Q)/E(\Q)_{\rm tors}$ up to finite index: replacing the basis by one of a sublattice of index $n$ multiplies both $\Reg_\infty$ and $\Reg_p$ by $n^2$, which leaves $q_\infty d_p$ unchanged. Reading the constant as $C\cdot\#\Sha$ requires that the basis generate $E(\Q)/E(\Q)_{\rm tors}$.
\end{remark}

\subsection{Results for \texorpdfstring{$y^2=x^3+39x$}{y2=x3+39x}}\label{subsec:obs39}
Here $C=8/4=2$ and $q_\infty=2$ numerically. The values $d_p$ were computed in two ways: by a program working with division polynomials modulo $p^2$ (no PARI), and with PARI's \texttt{ellpadicheight} and \texttt{ellpadicregulator}; the residues agree at all $1632$ primes. At $p=251,1091,2659$ SageMath's \texttt{padic\_regulator}, which does not use PARI, gives the same $d_p$.

\begin{observation}[numerical evidence]\label{obs:39}
For all $1632$ inert primes $7\le p<30000$, $d_p\not\equiv0\pmod p$ and $\Aplus\equiv-2d_p\pmod p$. At $p=7$ and $p=11$, $a_+\equiv13\pmod{49}$ and $a_+\equiv19\pmod{121}$, as predicted by the second congruence in \eqref{eq:prediction} with $(u,e',\ell')\equiv(25,25,22)\pmod{49}$ and $(106,96,56)\pmod{121}$.
\end{observation}

The residues of $a_+$ modulo $p^2$ at $p=7,11$ were obtained from the third layer of the cyclotomic tower, i.e.\ from the plus modular symbols $\msp{a/p^4}$ (the same eclib object) and the interpolation formulas \cite[(10)--(11)]{PollackRubin} at the characters of $\Gamma$ of order dividing $p^3$; we do not write out this bookkeeping, and the second-digit statements are numerical, like the rest of this section. The first values are $(d_p,\Aplus)=(4,6),(7,8),(5,9),(4,15),(28,6),(8,27),(31,32)$ for $p=7,11,19,23,31,43,47$.

Two remarks on the status of Observation~\ref{obs:39}. First, the congruence was noticed in the data before Proposition~\ref{prop:derivation} was written down; the $1632$ values were not predictions fixed in advance. Second, the two sides are computed from disjoint inputs (modular symbols on one side, the points $P,Q$ and $p$-adic heights on the other); they share only the model of $E$ and the choice of $\gamma$. If $\Aplus$ were wrong at some prime and $d_p$ were not, the congruence would hold there with probability about $1/p$.

As a check of the whole chain of normalisations, we compared both coordinates of \eqref{eq:BPR} on $E$ at $p=7,11,19,23$ with SageMath's implementation of the $\Dp$-valued $p$-adic $L$-series and regulator (\texttt{Dp\_valued\_series}, \texttt{padic\_regulator}): with $C\cdot\#\Sha=2$, the ratios of the two sides are $1+O(p)$ in both coordinates. This check is independent of our programs, but both sides come from the same software, and its precision is only $O(p)$.

\section{Two control curves}\label{sec:control}

\subsection{Choice of curves and procedure}
To test whether the constant in Observation~\ref{obs:39} really follows $q_\infty$, we used two curves of the same family with different predicted constants:
\[
E_{25721}\colon y^2=x^3+25721x,\qquad E_{51137}\colon y^2=x^3+51137x .
\]
Their basic data are in Table~\ref{tab:control-data}. In the notation $y^2=x^3-Dx$ these are $D=-25721$ and $D=-51137$. The curve $E_{25721}$ appears in \cite[Table~2]{CalongeFabregat}, where the identity of \cite[Thm.~4.3]{CalongeFabregat} is checked for it at the split prime $p=5$ to two $5$-adic digits; $E_{51137}$ is, according to the data archived with \cite{CalongeFabregat}, one of the $20$ rank-two curves with analytic $\Sha$ of order $9$ considered in \cite[\S8.2]{CalongeFabregat}. Here we add tests at inert primes.

\begin{table}[ht]
\centering\small
\renewcommand{\arraystretch}{1.2}
\begin{tabular}{@{}lll@{}}
\toprule
& $E_{25721}$ & $E_{51137}$\\
\midrule
conductor & $146506816=2^6\cdot17^2\cdot89^2$ & $167359537216=2^6\cdot51137^2$\\
minimal discriminant & $-2^6\cdot17^6\cdot89^3$ & $-2^6\cdot51137^3$\\
bad primes: type, $c_\ell$ & $2$: II, $1$;\ $17$: I$_0^*$, $4$;\ $89$: III, $2$ & $2$: II, $1$;\ $51137$: III, $2$\\
$E(\Q)_{\rm tors}$; root number & $\Z/2$; $+1$ & $\Z/2$; $+1$\\
$C=\prod c_\ell/\#E(\Q)_{\rm tors}^2$ & $2$ & $1/2$\\
rank & $2$ & $2$\\
basis $P$ & $(833/4,\ 30345/8)$ & $(256/225,\ 814096/3375)$\\
basis $Q$ & $(2448,\ -121380)$ & $\bigl(\tfrac{26615281}{115600},\ -\tfrac{192461494071}{39304000}\bigr)$\\
$\Omega_E$ & $0.2928095625196648089\ldots$ & $0.2465888337131131514\ldots$\\
$\Reg_\infty$ & $17.680168473040045367\ldots$ & $127.09314935370395738\ldots$\\
$L''(E,1)$ & $82.830758333837213834\ldots$ & $282.05776324850723403\ldots$\\
$q_\infty=L''(E,1)/(2\Omega_E\Reg_\infty)$ & $8.000\ldots$ & $4.500\ldots$\\
analytic order of $\Sha$, $q_\infty/C$ & $4$ & $9$\\
\bottomrule
\end{tabular}
\caption{The control curves. $q_\infty$ was computed to $60$ digits; for $E_{51137}$ an interval computation gives $|q_\infty-9/2|<1.2\cdot10^{-68}$.}\label{tab:control-data}
\end{table}

For both curves the rank is $2$ by a $2$-descent (for $E_{25721}$ a second descent is needed), the given points are saturated, and the model is minimal. A $2$-descent gives $\Sha(E_{25721}/\Q)[2]\supseteq(\Z/2)^2$ (cf.\ \cite[\S8.2]{CalongeFabregat} for curves of analytic order $4$), whereas the $2$-Selmer group of $E_{51137}$ has dimension $3$, so $\Sha(E_{51137}/\Q)[2]=0$; for the $20$ curves of analytic order $9$ in its range, which by the archived data include $E_{51137}$, \cite[\S8.2]{CalongeFabregat} reports an unconditional $3$-descent giving $\Sha(E/\Q)[3]\supseteq(\Z/3)^2$. The analytic orders $4$ and $9$ are numerical identifications.

\emph{Procedure.} The predictions were fixed in six registrations, each made before any value of $\Aplus$ in the corresponding batch was computed (for $E_{25721}$: $p=11$; then $p=7,19,23,31,43,47$; then the second digit at $p=11$; then the second digit at $p=7$; for $E_{51137}$: $p=7,11,19$; then $p=23,31,43$). Each registration recorded, with a time stamp: the curve, the model, the differential and period, the basis, the choice $\gamma\leftrightarrow1+p$, the sign convention; the values $d_p$ (and, for the second digit, $u\bmod p^2$, $e'\bmod p$), computed from the $p$-adic regulator; the predicted residues for the main constant $q_\infty$ and for alternative constants; and the rule that each prediction is compared separately and that no normalisation is adjusted after seeing a target value. Each later registration was made when the results of the earlier batches were already known. As a control, every analytic computation first reproduced the known values of $\Aplus$ for $y^2=x^3+39x$. The alternatives were: the constant $2$ of the original curve, and, for $E_{51137}$, the constant $1/2=C$ (i.e.\ discarding the factor $9$). In one registration the record stated that the heights had been obtained by two methods when at that moment only one had been run; the second method gave identical values before the corresponding analytic computation started.

The analytic side ($\Aplus$ on the control curves) was computed from the first and, for the second digit, the third layer of the cyclotomic tower. For $E_{25721}$ at all seven primes and for $E_{51137}$ at $p=7,11,19$, there were two interval-arithmetic implementations (with rational reconstruction based on an a priori bound for the denominators) and one independent floating-point implementation; for $E_{51137}$ at $p=23,31,43$ only one floating-point implementation without interval certification was available. The $d_p$ were computed with PARI and, independently, with a second implementation of $p$-adic heights.

\subsection{Results}

\begin{table}[ht]
\centering\small
\renewcommand{\arraystretch}{1.15}
\begin{tabular}{@{}rrrrrr@{}}
\toprule
$p$ & $d_p$ & $R_1/p^2\bmod p^2$ & predicted $-8d_p$ & alternative $-2d_p$ & observed $\Aplus\bmod p$\\
\midrule
7 & 0 & 28 & 0 & 0 & 0\\
11 & 2 & 90 & 6 & 7 & 6\\
19 & 12 & 31 & 18 & 14 & 18\\
23 & 15 & 475 & 18 & 16 & 18\\
31 & 24 & 272 & 25 & 14 & 25\\
43 & 5 & 1682 & 3 & 33 & 3\\
47 & 9 & 949 & 22 & 29 & 22\\
\bottomrule
\end{tabular}
\caption{$E_{25721}$, $q_\infty\approx8$: the first seven good supersingular primes $p\ge7$.}\label{tab:25721}
\end{table}

At the second digit: for $p=11$ the predictions were $105$ (constant $8$) and $117$ (constant $2$) modulo $121$, and $a_+\equiv105\pmod{121}$ was observed; for $p=7$, where the first digit does not separate the constants, the predictions were $28$ and $7$ modulo $49$, and $a_+\equiv28\pmod{49}$ was observed.

\begin{table}[ht]
\centering\small
\renewcommand{\arraystretch}{1.15}
\begin{tabular}{@{}rrrrrrr@{}}
\toprule
$p$ & $d_p$ & $R_1/p^2\bmod p^2$ & predicted $-\tfrac92d_p$ & alternative $-2d_p$ & alternative $-\tfrac12d_p$ & observed\\
\midrule
7 & 2 & 23 & 5 & 3 & 6 & 5\\
11 & 10 & 87 & 10 & 2 & 6 & 10\\
19 & 13 & 89 & 8 & 12 & 3 & 8\\
23 & 1 & 24 & 7 & 21 & 11 & 7$^\dagger$\\
31 & 9 & 939 & 6 & 13 & 11 & 6$^\dagger$\\
43 & 29 & 1749 & 20 & 28 & 7 & 20$^\dagger$\\
\bottomrule
\end{tabular}
\caption{$E_{51137}$, $q_\infty\approx9/2$: the first six good supersingular primes $p\ge7$. $^\dagger$\,One floating-point implementation only.}\label{tab:51137}
\end{table}

\subsection{What the tests show}
If $\Aplus\equiv-k\,d_p\pmod p$ with a rational constant $k$ independent of $p$, the data for $E_{25721}$ (the first digit at $p=19,23,31,43,47$ and both digits at $p=11$; the prime $7$, where $d_7=0$, is not used) force $k\equiv8\pmod{3312796927}$, where $3312796927=11^2\cdot19\cdot23\cdot31\cdot43\cdot47$; among rationals $a/b$ in lowest terms with $|a|,b\le40697$ (so $|a|,b<\sqrt{3312796927/2}$) the only solution is $k=8$. For $E_{51137}$ the first digit at the six primes forces $k\equiv9/2\pmod{44854117}$, $44854117=7\cdot11\cdot19\cdot23\cdot31\cdot43$; among rationals $a/b$ in lowest terms with $|a|,b\le4077646$ the only one in this residue class is $9/2$ (the next one has $b=4077649$). The second bound is far beyond $\sqrt{44854117/2}$, so this is not a standard rational reconstruction; it was checked by exhaustive search over $b$. (Note $9/2\equiv8\equiv1\pmod7$, so the prime $7$ alone does not separate these constants.) No single constant works for all three curves: at $p=11$, $E_{25721}$ requires $k\equiv8$ and $E_{51137}$ requires $k\equiv10\pmod{11}$.

Thus, at the first $p$-adic digit (and at the second digit in two cases), the observed coefficients are consistent with the constant $q_\infty$ predicted by Proposition~\ref{prop:derivation} combined with the complex BSD formula, and inconsistent with the alternatives fixed in advance. In particular, for $E_{51137}$ the constant $C=1/2$, which would ignore the odd factor $9$ of the analytic order of $\Sha$, fails at all six primes.

\emph{What the tests do not show.} They do not measure $\#\Sha$: they compare the $p$-adic leading coefficient with the complex BSD quotient. If one takes the analytic order of $\Sha$ as input, the data test the $p$-adic formula; if one wants to bound $\#\Sha$ from the data, one must assume the $p$-adic formula. These two readings cannot serve as independent evidence for each other. For example, at $p=11$ on $E_{25721}$, the observed residue $6$ only requires $\#\Sha\equiv4\pmod{11}$, given the $p$-adic formula and finiteness of $\Sha$. Each test is a single $p$-adic digit (two digits in four cases: $p=7,11$ for $y^2=x^3+39x$ and for $E_{25721}$), on three curves.

\subsection{Relation to previous work}
Calonge Fabregat \cite{CalongeFabregat} formulates and tests an exact $p$-adic BSD quotient for rank-two curves $y^2=x^3-Dx$ at split primes $p\equiv1\pmod4$, $p\le41$: in $5520$ pairs $(E,p)$ with conductors up to $3.9\cdot10^{12}$ the agreement holds modulo $p^3$, including curves with analytic $\Sha$ of order $4$ and $9$ \cite[\S1.3, \S8.2, Table~1]{CalongeFabregat}; a further identity, which concerns the position of Rubin's elliptic-unit class and involves $q_\infty$ \cite[Thm.~4.3]{CalongeFabregat}, shows no discrepancy at $p=5$ for $301$ curves with analytic $\Sha$ of order $4$ or $9$ and $q_\infty\in\{4,8,9,16,18,32,9/2\}$ \cite[\S1.3, \S8.4]{CalongeFabregat}. That work treats only split primes; its Table~2 includes our curve $E$ at $p=5$, where the identity of \cite[Thm.~4.3]{CalongeFabregat} holds to five $5$-adic digits. What we add is the corresponding test at inert primes, at the first digit only and on three curves. The curve $E$ itself, with the same generators $P,Q$ and their heights, already appears in the tables of Bernardi, Goldstein and Stephens \cite[p.~168]{BGS84}; we have read only that page of \cite{BGS84}.

At supersingular primes, absolute leading-term comparisons with the $p$-adic regulator have been made before for curves without complex multiplication: by Perrin-Riou \cite[\S6.4 of the arXiv version]{PR03} at $p=3$ for several curves of rank two, and by Stein and Wuthrich \cite[\S12.2, p.~1786]{SteinWuthrich} for the curve \texttt{1483a1} at $p=5$. In both works the comparison, together with Kato's divisibility, was also used to determine the $p$-part of $\Sha$; for \texttt{1483a1}, Stein and Wuthrich conclude that $\Sha(E/\Q)[5]=0$. Kurihara and Pollack \cite[\S0.1]{KuriharaPollack} computed, for quadratic twists of $X_0(17)$ and $X_0(32)$ at $p=3$, the ratio of the leading coefficients of the two signed $p$-adic $L$-functions and compared it with $(p-1)/2$ (a criterion for the Selmer rank). Pollack \cite[\S7, Table~1]{Pollack} tabulated the signed analytic Iwasawa invariants of quadratic twists of $X_0(32)$, which have complex multiplication by $\Z[i]$, at the inert prime $p=3$, including twists of rank two. On our reading of that table, for some of the rank-two twists (for example the twist by $65$) the plus $\lambda$-invariant is $2$ with both zeros at $T=0$; if moreover $\mu^+=0$, which the table does not list explicitly, the argument of Proposition~\ref{prop:criterion}, combined with \cite[Thm.~7.3]{PollackRubin} (which appeared later), would apply to such a twist at $p=3$, subject to the local hypotheses at its bad primes. We have not checked this; it is not stated in \cite{Pollack}. Neither \cite{KuriharaPollack} nor \cite{Pollack} compares a rank-two leading coefficient with the regulator (in rank one, Kurihara and Pollack \cite[\S0.2]{KuriharaPollack} use the $p$-adic BSD conjecture to construct rational points from the derivative). We have not been able to read \cite{BPR93}, which according to its review contains numerical evidence for one curve with supersingular reduction.

\section{Scope and limitations}\label{sec:limits}

\begin{enumerate}
\item \emph{Finitely many primes.} Theorems~A and~B concern primes below $30000$. They imply nothing about the finiteness of $\Sha(E/\Q)$ or about the BSD formula for $E$.
\item \emph{Computer assistance.} Theorem~A depends on the correctness of the computation of Section~\ref{sec:computation}, and in particular of eclib's modular symbols of level $48672$, which are the single input to all $1632$ values (Remark~\ref{rem:shared}). Propositions~\ref{prop:p2}, \ref{prop:p3}, \ref{prop:p13} depend on mwrank, on PARI's \texttt{bnfcertify}, and on an exact finite-field computation, respectively.
\item \emph{Inputs from the literature.} Theorem~B uses \cite[Thm.~1.3]{CLS2}, whose computation was done in MAGMA and is only partly tabulated, and the preprint \cite{Banwait} at $p=17$. The criterion uses the theorems of \cite{PollackRubin}, which in turn rest on \cite{Kobayashi,Pollack,Rubin91} and further work cited there.
\item \emph{Only the plus sign.} The criterion uses $L_p^+$ only. We make no claim about $L_p^-$.
\item \emph{The criterion is one-sided.} $\Aplus\in\Zp^\times$ is equivalent to $X^+\cong\Zp^2$, which is stronger than $\Sha(E/\Q)[p^\infty]=0$; a prime with $\Aplus\equiv0\pmod p$ would not be a counterexample to anything. There is no such prime below $30000$.
\item \emph{Numerical identifications.} The values $q_\infty=2,8,9/2$ (and the analytic orders $1,4,9$) are numerical. For $E_{51137}$ there is an interval enclosure; for the other two curves only high-precision floating-point values.
\item \emph{Part II is evidence, not proof.} Observation~C and Section~\ref{sec:control} are consistent with Conjecture~\ref{conj:BPR}, at one $p$-adic digit per prime (two digits in four cases). For $y^2=x^3+39x$ the $1632$ congruences were observed, not predicted in advance; the control-curve predictions were fixed in advance. For $E_{51137}$ at $p=23,31,43$ the analytic values come from a single floating-point implementation.
\item \emph{Normalisation of $d_p$.} We assume that PARI's and SageMath's $p$-adic regulators agree with \cite[Lemma~4.2]{SteinWuthrich}; the agreement of both coordinates in the check at the end of Section~\ref{sec:leading} supports this.
\end{enumerate}

\section*{Data and code}

The ancillary table \texttt{inert\_primes.csv} lists $\Aplus\bmod p$, $d_p\bmod p$, $-2d_p\bmod p$, and the outcome of the comparison $\Aplus\equiv-2d_p\pmod p$ for each of the $1632$ good supersingular primes $7\le p<30000$. The accompanying \texttt{README.md} defines the columns and normalisations. In \texttt{code/}, \texttt{build\_edges.py} constructs the exact modular-symbol table of Section~\ref{subsec:exact} using SageMath~10.9 \cite{sagemath} and eclib~20250627 \cite{eclib}. The C++ program \texttt{plus\_coefficient.cpp} uses this table to compute $\Aplus\bmod p$ by the finite sum in Section~\ref{sec:computation}. The Python program \texttt{height\_digit.py} computes $d_p$ by division polynomials as in Section~\ref{subsec:obs39}. The program \texttt{sample/run.py} rebuilds the symbols without eclib and repeats the $20$-prime check in Section~\ref{subsec:checks}, item~4(d). The programs \texttt{control\_25721.c} and \texttt{control\_51137.c} reproduce the leading-term comparisons for the respective curves in Section~\ref{sec:control}. The file \texttt{code/README.md} gives the dependencies, commands, expected output and approximate running times.

\section*{Use of generative AI}

In accordance with the arXiv policy on the use of generative AI tools, we report that such tools were used substantially in this work. The computer programs, the computations, the consistency checks between independent implementations, the comparison with the literature, and a first complete draft of the text, including the written proofs, were produced by AI systems based on large language models, working under the authors' direction. Two independent AI systems were used to cross-check each other's computations and arguments. The authors take full responsibility for the content of this paper.

\bibliographystyle{amsplain}
\bibliography{refs,extra}

\providecommand{\url}[1]{\texttt{#1}}
\providecommand{\bysame}{\leavevmode\hbox to3em{\hrulefill}\thinspace}
\providecommand{\MR}{\relax\ifhmode\unskip\space\fi MR }
\providecommand{\MRhref}[2]{%
  \href{http://www.ams.org/mathscinet-getitem?mr=#1}{#2}
}
\providecommand{\href}[2]{#2}
\begin{thebibliography}{10}

\bibitem{Banwait}
Barinder~S. Banwait, \emph{Second derivatives of $p$-adic {$L$}-functions and
  the {S}hafarevich--{T}ate group of rank-two {CM} elliptic curves}, 2026,
  Preprint, \texttt{arXiv:2609.08431v1} [math.NT].

\bibitem{BGS84}
Dominique Bernardi, Catherine Goldstein, and Nelson Stephens, \emph{Notes
  $p$-adiques sur les courbes elliptiques}, J. Reine Angew. Math. \textbf{351}
  (1984), 129--170.

\bibitem{BPR93}
Dominique Bernardi and Bernadette Perrin-Riou, \emph{Variante $p$-adique de la
  conjecture de {B}irch et {S}winnerton-{D}yer (le cas supersingulier)}, C. R.
  Acad. Sci. Paris S\'er. I Math. \textbf{317} (1993), no.~3, 227--232.

\bibitem{CalongeFabregat}
Miquel Calonge~Fabregat, \emph{Elliptic units, derived {H}eegner classes and
  the $p$-adic {B}irch and {S}winnerton-{D}yer quotient of {CM} elliptic curves
  of rank two}, 2026, Preprint, \texttt{arXiv:2611.09192v1} [math.NT]. Code and
  data: \url{https://doi.org/10.5281/zenodo.23982770}.

\bibitem{Castella25}
Francesc Castella, \emph{A formula of {P}errin-{R}iou and characteristic power
  series of signed {S}elmer groups}, 2025, Preprint,
  \texttt{arXiv:2502.19618v2} [math.NT].

\bibitem{CLS1}
John Coates, Zhibin Liang, and Ramdorai Sujatha, \emph{The
  {T}ate--{S}hafarevich group for elliptic curves with complex multiplication},
  J. Algebra \textbf{322} (2009), no.~3, 657--674.

\bibitem{CLS2}
\bysame, \emph{The {T}ate--{S}hafarevich group for elliptic curves with complex
  multiplication {II}}, Milan J. Math. \textbf{78} (2010), no.~2, 395--416.

\bibitem{Cremona97}
J.~E. Cremona, \emph{Algorithms for modular elliptic curves}, second ed.,
  Cambridge University Press, Cambridge, 1997.

\bibitem{eclib}
\bysame, \emph{eclib}, \url{https://github.com/JohnCremona/eclib}, 2025,
  Version 20250627. DOI 10.5281/zenodo.593635.

\bibitem{deShalit}
Ehud de~Shalit, \emph{Iwasawa theory of elliptic curves with complex
  multiplication}, Perspectives in Mathematics, vol.~3, Academic Press, Boston,
  MA, 1987.

\bibitem{Kobayashi}
Shin-ichi Kobayashi, \emph{Iwasawa theory for elliptic curves at supersingular
  primes}, Invent. Math. \textbf{152} (2003), no.~1, 1--36.

\bibitem{KuriharaPollack}
Masato Kurihara and Robert Pollack, \emph{Two $p$-adic {$L$}-functions and
  rational points on elliptic curves with supersingular reduction},
  {$L$}-functions and {G}alois representations (David Burns, Kevin Buzzard, and
  Jan Nekov\'a\v{r}, eds.), London Math. Soc. Lecture Note Ser., vol. 320,
  Cambridge Univ. Press, Cambridge, 2007, pp.~300--332.

\bibitem{PR03}
Bernadette Perrin-Riou, \emph{Arithm\'etique des courbes elliptiques \`a
  r\'eduction supersinguli\`ere en $p$}, Experiment. Math. \textbf{12} (2003),
  no.~2, 155--186.

\bibitem{Pollack}
Robert Pollack, \emph{On the $p$-adic {$L$}-function of a modular form at a
  supersingular prime}, Duke Math. J. \textbf{118} (2003), no.~3, 523--558.

\bibitem{PollackRubin}
Robert Pollack and Karl Rubin, \emph{The main conjecture for {CM} elliptic
  curves at supersingular primes}, Ann. of Math. (2) \textbf{159} (2004),
  no.~1, 447--464.

\bibitem{Rubin91}
Karl Rubin, \emph{The ``main conjectures'' of {I}wasawa theory for imaginary
  quadratic fields}, Invent. Math. \textbf{103} (1991), no.~1, 25--68.

\bibitem{Silverman}
Joseph~H. Silverman, \emph{The arithmetic of elliptic curves}, second ed.,
  Graduate Texts in Mathematics, vol. 106, Springer, New York, 2009.

\bibitem{Sprung}
Florian Sprung, \emph{A formulation of $p$-adic versions of the {B}irch and
  {S}winnerton-{D}yer conjectures in the supersingular case}, Res. Number
  Theory \textbf{1} (2015), Paper No. 17, 13 pp.

\bibitem{SteinWuthrich}
William Stein and Christian Wuthrich, \emph{Algorithms for the arithmetic of
  elliptic curves using {I}wasawa theory}, Math. Comp. \textbf{82} (2013),
  no.~283, 1757--1792.

\bibitem{PARI2}
{The PARI~Group}, Univ. Bordeaux, \emph{{PARI/GP version \texttt{2.17.3}}},
  2025, available from \url{http://pari.math.u-bordeaux.fr/}.

\bibitem{sagemath}
{The Sage Developers}, \emph{{S}age{M}ath, the {S}age {M}athematics {S}oftware
  {S}ystem ({V}ersion 10.9)}, 2026, DOI 10.5281/zenodo.8042260.

\end{thebibliography}

\end{document}